\documentclass{amsart}
\pdfoutput=1

\usepackage[margin=1in]{geometry}
\usepackage[utf8]{inputenc}
\usepackage{graphicx,cite}
\usepackage{amsmath,amssymb,amsthm}
\usepackage{color}
\usepackage{tikz}
\usetikzlibrary{arrows.meta}
\usepackage{pgfplots}
\pgfplotsset{compat=1.18}
\usepgfplotslibrary{fillbetween}
\allowdisplaybreaks

\makeatletter
\def\l@section{\@tocline{1}{12pt plus2pt}{0pt}{}{\bfseries}}
\def\l@subsection{\@tocline{2}{0pt}{2pc}{2pc}{}}
\makeatother

\makeatletter
\def\subsection{\@startsection{subsection}{2}{\z@}%
  {-3.25ex\@plus -1ex \@minus -.2ex}%
  {1.5ex \@plus .2ex}%
  {\normalfont\bfseries\boldmath}}
\def\subsubsection{\@startsection{subsubsection}{3}%
  \z@{.5\linespacing\@plus.7\linespacing}{-.5em}%
  {\normalfont\bfseries\boldmath}}
\renewcommand\paragraph{\@startsection{paragraph}{4}{\z@}%
  {3.25ex \@plus1ex \@minus.2ex}%
  {-1em}%
  {\normalfont\normalsize\bfseries}}
\makeatother

\theoremstyle{plain}
\newtheorem{thm}{Theorem}[section]

\newtheorem{lem}[thm]{Lemma}
\newtheorem{prop}[thm]{Proposition}

\theoremstyle{definition}
\newtheorem{defn}[thm]{Definition}

\theoremstyle{remark}
\newtheorem{rem}[thm]{Remark}

\theoremstyle{plain}

\newtheorem{opques}[thm]{Open Question}

\numberwithin{equation}{section}

\newcommand{\thistheoremname}{}
\newtheorem{genericthm}[thm]{\thistheoremname}

\newtheorem*{genericthm*}{\thistheoremname}
\newenvironment{namedthm*}[1]
  {\renewcommand{\thistheoremname}{#1}%
   \begin{genericthm*}}
  {\end{genericthm*}}

\newcommand{\norm}[1]{\left\Vert#1\right\Vert}

\newcommand{\eps}{\varepsilon}

\newcommand{\R}{\mathbb R}

\newcommand{\C}{\mathbb C}

\newcommand{\N}{\mathbb N}

\newcommand{\Z}{\mathbb Z}

\newcommand{\calD}{\mathcal D}

\newcommand{\dist}{\hbox{\rm dist}}

\newcommand{\calS}{\mathcal S}
\newcommand{\supp}{\textnormal{supp}}

\makeatletter
\newcommand{\vast}{\bBigg@{4}}
\newcommand{\Vast}{\bBigg@{5}}
\makeatother

\def\udot#1{\ifmmode
  \oalign{$#1$\crcr\hidewidth.\hidewidth}
  \else
  \oalign{#1\crcr\hidewidth.\hidewidth}
  \fi}

\def\beq{\begin{equation}}
\def\eeq{\end{equation}}

\makeatletter
\newcommand{\doublewidetilde}[1]{{%
  \mathpalette\double@widetilde{#1}%
}}
\newcommand{\double@widetilde}[2]{%
  \sbox\z@{$\m@th#1\widetilde{#2}$}%
  \ht\z@=.9\ht\z@
  \widetilde{\box\z@}%
}
\makeatother

\def\one{\mbox{1\hspace{-4.25pt}%
  \fontsize{12}{14.4}\selectfont\textrm{1}}}

\makeatletter
\def\@makefnmark{%
  \leavevmode
  \raise.9ex\hbox{%
    \fontsize\sf@size\z@\normalfont\tiny\@thefnmark}}
\makeatother

\usepackage[hidelinks]{hyperref}

\title[Weak $L^{2/3}$ bounds for the bilinear Hilbert transform]
{Weak $L^{2/3}$ bounds for the bilinear Hilbert transform}
\author{Bingyang Hu}
\address{Department of Mathematics and Statistics\\Auburn University\\Auburn, Alabama, U.S.A., 36849}
\email{bzh0108@auburn.edu}
\date{November 9, 2026}
\begin{document}
\begin{abstract}
For Schwartz functions $f_1,f_2$ on $\R$, define the
bilinear Hilbert transform by
\[
 \textnormal{\bf BHT}(f_1,f_2)(x)
 =\operatorname{p.v.}\int_{\R}
 f_1(x-t)f_2(x+t)\,\frac{dt}{t}.
\]
We prove the weak-type estimate
\[
 \norm{\textnormal{\bf BHT}(f_1,f_2)}_{L^{2/3,\infty}}
 \lesssim_{p_1,p_2}\norm{f_1}_{p_1}\norm{f_2}_{p_2}
\]
for $1<p_1,p_2<2$ with $1/p_1+1/p_2=3/2$.
The key ingredient is a new local supremum estimate for packet sums, which we use to establish the forest estimate.
\end{abstract}
\maketitle
\tableofcontents

\section{Introduction}
For $f_1,f_2\in\calS(\R)$, we consider the bilinear Hilbert transform
\begin{equation}\label{eq:bht}
 \textnormal{\bf BHT}(f_1,f_2)(x)
 =\operatorname{p.v.}\int_{\R}
 f_1(x-t)f_2(x+t)\,\frac{dt}{t}.
\end{equation}
Our main result is the following.

\begin{thm}\label{thm:main}
Let
\begin{equation}\label{eq:critical}
 1<p_1,p_2<2,\qquad \textrm{with} \qquad \frac1{p_1}+\frac1{p_2}=\frac32.
\end{equation}
Then for any $f_1, f_2 \in \calS(\R)$, one has
\begin{equation}\label{eq:main}
 \norm{\textnormal{\bf BHT}(f_1,f_2)}_{L^{2/3,\infty}}
 \lesssim_{p_1,p_2}\norm{f_1}_{p_1}\norm{f_2}_{p_2}.
\end{equation}
Consequently, $\textnormal{\bf BHT}$ extends uniquely to a bounded bilinear operator
\[
 \textnormal{\bf BHT}:L^{p_1}(\R)\times L^{p_2}(\R)
 \longrightarrow L^{2/3,\infty}(\R).
\]
\end{thm}

The bilinear Hilbert transform arose from Calder\'on's celebrated
work on the Cauchy integral on Lipschitz curves~\cite{Calderon}.
Besides the translation and dilation symmetries familiar from
classical Calder\'on--Zygmund theory, $\textnormal{\bf BHT}$
also has \emph{modulation invariance}, a property shared
with the Carleson operator.
To study the resonance associated with modulation invariance,
one uses wave-packet analysis (i.e. time-frequency analysis) to keep track of the localization and interaction of the input functions in both space and frequency.

Building on the time--frequency ideas of
Carleson~\cite{Carleson} and Fefferman~\cite{Fefferman}
in their work on the pointwise convergence of Fourier series,
Lacey and Thiele established, in their breakthrough
papers~\cite{LT97,LT99}, the strong-type estimates
\begin{equation}\label{eq:lacey-thiele}
 \norm{\textnormal{\bf BHT}(f_1,f_2)}_r
 \lesssim_{q_1,q_2}
 \norm{f_1}_{q_1}\norm{f_2}_{q_2},
 \qquad
 \frac1r=\frac1{q_1}+\frac1{q_2},
\end{equation}
for $1<q_1,q_2<\infty$ and $2/3<r<\infty$.
In their proofs, Lacey and Thiele decompose the operator
into wave packets, organize the tiles into trees, and
combine estimates on individual trees with almost
orthogonality. Muscalu, Tao, and Thiele~\cite{MTT} developed these ideas
further to treat a class of multilinear operators whose
multipliers are singular along nondegenerate subspaces.
For a systematic account of wave-packet analysis,
including its applications to the Carleson operator and
the bilinear Hilbert transform, we refer to~\cite{Thiele}.

\vspace{0.1cm}

However, whether the strong-type estimate \eqref{eq:lacey-thiele}
extends to $1/2<r\le2/3$ remains a longstanding open
problem in harmonic analysis.
We formulate the question below;
see~\cite{DPT} and~\cite[Remark~1.9]{HL}.

\begin{opques}\label{ques:bht-strong}
Let $1<q_1,q_2<\infty$ and $1/2<r\le2/3$ satisfy
\[
 \frac1r=\frac1{q_1}+\frac1{q_2}.
\]
Does the strong-type estimate
\[
 \norm{\textnormal{\bf BHT}(f_1,f_2)}_r
 \lesssim_{q_1,q_2}
 \norm{f_1}_{q_1}\norm{f_2}_{q_2}
\]
hold for all $f_1,f_2\in\calS(\R)$?
\end{opques}

At the critical exponent $r=2/3$, partial progress toward
Open Question~\ref{ques:bht-strong} has been made.

\begin{itemize}
\item
Lacey~\cite{Lacey2000} showed that
uniform $L^{q_1}\times L^{q_2}\to L^r$ estimates for
bilinear wave-packet models with arbitrary signs
fail when $r<2/3$, where $q_1,q_2>1$ and
$1/q_1+1/q_2=1/r$.
His counterexample uses packets at a single spatial scale.

\item
Demeter and Di Plinio~\cite{DDP} then studied the
quartile operator, a Walsh model of $\textnormal{\bf BHT}$.
They proved weak-type estimates at $r=2/3$ when one input is bounded by an indicator function.

\item
Later, Di Plinio and Thiele~\cite{DPT} obtained endpoint
estimates for $\textnormal{\bf BHT}$ with logarithmic
corrections.
In particular, they proved generalized restricted weak-type
estimates with a doubly logarithmic loss.

\item
Subsequently, Culiuc, Di Plinio, and Ou~\cite{CDO}
adapted Lacey's construction to the class of bilinear
multipliers studied in~\cite{MTT}.
They constructed a family of symbols satisfying uniform
derivative bounds for which the
$L^{q_1}\times L^{q_2}\to L^r$ operator norms are unbounded
when $r<2/3$, with $q_1,q_2>1$ and
$1/q_1+1/q_2=1/r$.
\end{itemize}
It is worthwhile mentioning that the counterexamples in~\cite{Lacey2000,CDO}
do not rule out strong-type estimates for
$\textnormal{\bf BHT}$ itself in the range $1/2<r<2/3$.

In the present paper, we take an important step toward
resolving Open Question~\ref{ques:bht-strong} by proving
the weak-type bound at $r=2/3$ for every exponent pair
satisfying \eqref{eq:critical}, as stated in
Theorem~\ref{thm:main}.

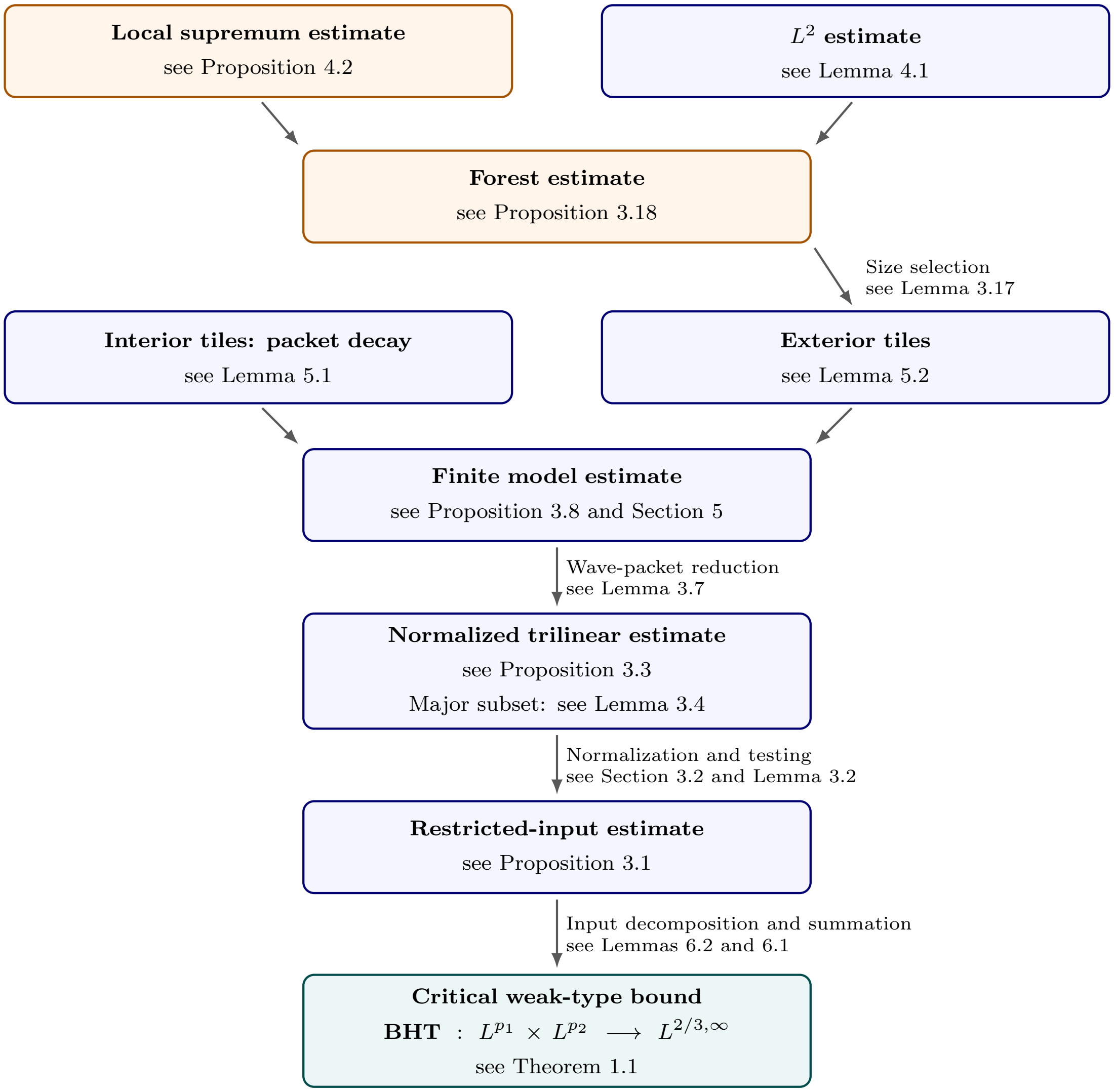
\begin{figure}[htbp]
\centering
\begin{tikzpicture}[
    font=\footnotesize,
    box/.style={
        rectangle,
        rounded corners=4pt,
        thick,
        draw=blue!45!black,
        fill=blue!4,
        align=center,
        text width=0.38\textwidth,
        minimum height=12mm,
        inner sep=5pt
    },
    core/.style={
        box,
        draw=orange!65!black,
        fill=orange!8
    },
    result/.style={
        box,
        draw=teal!60!black,
        fill=teal!7
    },
    arrow/.style={
        -{Latex[length=2mm]},
        thick,
        draw=black!65,
        shorten <=2pt,
        shorten >=2pt
    },
    explanation/.style={
        font=\scriptsize,
        align=left,
        inner sep=3pt
    }
]

\node[core] (local) at (-3.9,0) {
    \textbf{Local supremum estimate}\\[1mm]
    see Proposition~\ref{prop:local}
};
\node[box] (l2) at (3.9,0) {
    \textbf{$L^2$ estimate}\\[1mm]
    see Lemma~\ref{lem:forest2}
};

\node[core] (forest) at (0,-1.9) {
    \textbf{Forest estimate}\\[1mm]
    see Proposition~\ref{prop:forest-p}
};
\draw[arrow] (local.south) -- (forest.north west);
\draw[arrow] (l2.south) -- (forest.north east);

\node[box] (inside) at (-3.9,-4.0) {
    \textbf{Interior tiles: packet decay}\\[1mm]
    see Lemma~\ref{lem:interior}
};
\node[box] (outside) at (3.9,-4.0) {
    \textbf{Exterior tiles}\\[1mm]
    see Lemma~\ref{lem:exterior}
};
\draw[arrow] (forest.south east)
    -- node[explanation,right,xshift=3mm] {
        Size selection\\
        see Lemma~\ref{lem:size-selection}
    } (outside.north);

\node[box] (model) at (0,-5.8) {
    \textbf{Finite model estimate}\\[1mm]
    see Proposition~\ref{prop:model} and Section~\ref{sec:model-proof}
};
\draw[arrow] (inside.south) -- (model.north west);
\draw[arrow] (outside.south) -- (model.north east);

\node[box] (normalized) at (0,-8.1) {
    \textbf{Normalized trilinear estimate}\\[1mm]
    see Proposition~\ref{prop:normalized}\\[1mm]
    Major subset: see Lemma~\ref{lem:major-subset}
};
\draw[arrow] (model.south)
    -- node[explanation,right] {
        Wave-packet reduction\\
        see Lemma~\ref{lem:model-reduction}
    } (normalized.north);

\node[box] (restricted) at (0,-10.4) {
    \textbf{Restricted-input estimate}\\[1mm]
    see Proposition~\ref{prop:restricted}
};
\draw[arrow] (normalized.south)
    -- node[explanation,right] {
        Normalization and testing\\
        see Section~\ref{subsec:normalization}
        and Lemma~\ref{lem:testing}
    } (restricted.north);

\node[result] (main) at (0,-12.8) {
    \textbf{Critical weak-type bound}\\[1mm]
    $\displaystyle
      \textnormal{\bf BHT}:L^{p_1}\times L^{p_2}
      \longrightarrow L^{2/3,\infty}
    $\\[1mm]
    see Theorem~\ref{thm:main}
};
\draw[arrow] (restricted.south)
    -- node[explanation,right] {
        Input decomposition and summation\\
        see Lemmas~\ref{lem:quantiles}
        and~\ref{lem:superposition}
    } (main.north);

\end{tikzpicture}
\caption{Structure of the proof of Theorem~\ref{thm:main}.}
\label{fig:proof-roadmap}
\end{figure}

The proof of Theorem~\ref{thm:main} uses the wave-packet
framework of~\cite{MTT} and follows the reductions and
size summation in~\cite{DDP,DPT}.
The main \emph{new} ingredient is the forest estimate in
Proposition~\ref{prop:forest-p}.
To prove the forest estimate, we combine the local supremum
estimate in Proposition~\ref{prop:local} with the $L^2$ bound
in Lemma~\ref{lem:forest2}. Since the proof involves several reductions, we summarize
its structure below (see Figure \ref{fig:proof-roadmap}).

\begin{rem}
The argument for the local supremum estimate in
Proposition~\ref{prop:local} can also be adapted to the
wave-packet models associated with the bilinear multiplier
class of Muscalu, Tao, and Thiele~\cite{MTT}.
Indeed, the proof uses only the tri-tile geometry and the
wave-packet bounds in
\eqref{eq:packet}--\eqref{eq:packet-deriv}, which are also
available in the discretization of this bilinear multiplier class;
see~\cite{MTT,CDO}.
For clarity, however, we keep the discussion focused on
$\textnormal{\bf BHT}$ and do not pursue the general
multiplier setting here.
\end{rem}

The rest of the paper is organized as follows.
In Section~\ref{sec:prelim}, we fix notation and recall
the necessary preliminaries.
In Section~\ref{sec:reduction}, we reduce
Proposition~\ref{prop:restricted} to a wave-packet model
and introduce trees, sizes, and the size selection lemma.
In Section~\ref{sec:local}, we prove the local supremum
estimate and then the forest estimate.
Section~\ref{sec:model-proof} contains the proof of
Proposition~\ref{prop:model}, which completes the proof
of Proposition~\ref{prop:restricted}.
Finally, in Section~\ref{sec:completion}, we remove the
restriction on the second input and complete the proof
of Theorem~\ref{thm:main}.

Throughout this paper, for any two nonnegative quantities $a$ and $b$, we write $a \lesssim b$ if there exists an absolute constant $C > 0$ independent of $a$ and $b$ such that $a \le Cb$. The notation $a \simeq b$ indicates that both $a \lesssim b$ and $b \lesssim a$ hold simultaneously.

\begin{rem}
While this manuscript was in preparation, we were pleased
to see a preprint by Di Plinio~\cite{DP2026} appear on arXiv,
establishing the same endpoint estimate as in
Theorem~\ref{thm:main}.
The present paper provides an alternative proof of the
endpoint estimate, based on a different treatment of the
key forest estimate.
Di Plinio~\cite{DP2026} proves the forest estimate through
an elegant refinement of the multi-frequency
Calder\'on--Zygmund decomposition, while our proof relies
on the local supremum estimate in
Proposition~\ref{prop:local}.
\end{rem}

\bigskip 

\section{Preliminaries}\label{sec:prelim}

\subsection{Notation and maximal functions}
For an interval $I$, let $c(I)$ denote its center and $|I|$ its length.
For $b>0$, let $bI$ be the interval with center $c(I)$ and length
$b|I|$. We write $\one_E$ for the characteristic function of a measurable set $E$,
and $p'=p/(p-1)$ for the conjugate exponent when $1<p<\infty$.
We also use the conventions
\[
 \left\langle f,g\right\rangle=\int_{\R}f(x)\overline{g(x)}\,dx,\qquad
 \widehat f(\xi)=\int_{\R}f(x)e^{-2\pi i x\xi}\,dx, 
\]
and denote the
inverse Fourier transform by
$\check u(x)=\int_\R u(\xi)e^{2\pi i x\xi}\,d\xi$. For a locally integrable function $f$, define
\[
 Mf(x)=\sup_{I\ni x}\frac1{|I|}\int_I|f(y)|\,dy,\qquad \textrm{and} \qquad 
 M_pf=(M(|f|^p))^{1/p}
\]
to be the standard Hardy--Littlewood maximal operator and its
$L^p$ variant, respectively.

\subsection{Dyadic lattices and tiles}
\begin{defn}[Dyadic lattice]\label{def:grid}
A \emph{dyadic lattice} on $\R$ is a collection
$\calD=\bigcup_{k\in\Z}\calD_k$ with the following properties:
\begin{enumerate}
\item[(i)] For some fixed $\ell>0$, $\calD_k$ is a partition of
$\R$ into half-open intervals of length $\ell 2^{-k}$.
\item[(ii)] Each interval in $\calD_k$ is the disjoint union
of two intervals in $\calD_{k+1}$.
We call these two intervals its \emph{children}.
\end{enumerate}
\end{defn}

\begin{defn}[Tiles and tri-tiles]\label{def:tile}
Let $I,\omega\subset\R$ be intervals satisfying
$|I||\omega|=1$.
We call the rectangle $I\times\omega$ a \emph{tile},
with spatial interval $I$ and frequency interval $\omega$.

A \emph{tri-tile} is a triple $s=(s_1,s_2,s_3)$ of tiles
with a common spatial interval $I_s$.
For each $i\in\{1,2,3\}$, we write $s_i=I_s\times\omega_{s,i}$. 
\end{defn}

\subsection{Wave packets}
To each tile $I\times\omega$, we associate a smooth
$L^1$-normalized \emph{wave packet} $\varphi$
whose Fourier transform is supported in $\omega$
and which decays away from $I$. More precisely,
we fix sufficiently large integers $A$ and $m_0$, and require that
\begin{equation}\label{eq:packet}
 \supp\,\widehat{\varphi}\subset\omega,\qquad
 |\varphi(x)|\lesssim |I|^{-1}
       \left(1+\frac{\dist(x,I)}{|I|}\right)^{-A},
\end{equation}
and
\begin{equation}\label{eq:packet-deriv}
 \left|\frac{d^m}{dx^m}
       \left(e^{-2\pi i c(\omega)x}\varphi(x)\right)\right|
 \lesssim |I|^{-1-m}
       \left(1+\frac{|x-c(I)|}{|I|}\right)^{-A},
 \qquad 0\le m\le m_0.
\end{equation}
The implicit constants in \eqref{eq:packet} and
\eqref{eq:packet-deriv} are fixed throughout the paper,
including in all suprema over wave-packet families. 

For example, we may take
\[
 \varphi(x)=|I|^{-1}e^{2\pi i c(\omega)x}
           \phi\!\left(\frac{x-c(I)}{|I|}\right),
\]
where $\phi$ is a fixed Schwartz function satisfying
$\supp\,\widehat{\phi}\subset[-1/2,1/2]$.
Observe that 
$\norm{\varphi}_1\lesssim1$ and $\norm{\varphi}_2\lesssim |I|^{-1/2}$.
For the tile $s_i$ of a tri-tile, we denote the corresponding wave packet by
$\varphi_{s,i}$.

To this end, observe that one can bound
$|\langle f,\varphi_{s,i}\rangle|$ by a constant multiple of
$\inf_{y\in I_s}Mf(y)$. Indeed, by decomposing $\R$ into
$2I_s$ and the annuli
$2^{m+1}I_s\setminus 2^m I_s$, $m\ge1$, we have
\begin{align}
 |\langle f,\varphi_{s,i}\rangle|
 &=\left|\int_{\R}f(x)\overline{\varphi_{s,i}(x)}\,dx\right|
 \notag\\
 &\lesssim \frac1{|I_s|}\int_{2I_s}|f(x)|\,dx
 +\sum_{m\ge1}\frac{2^{-mA}}{|I_s|}
   \int_{2^{m+1}I_s\setminus 2^m I_s}|f(x)|\,dx
 \notag\\
 &\lesssim \sum_{m\ge0}2^{-m(A-1)}
   \frac1{|2^{m+1}I_s|}
   \int_{2^{m+1}I_s}|f(x)|\,dx
 \notag\\
 &\lesssim \inf_{y\in I_s}Mf(y),
 \label{eq:packet-max}
\end{align}
which gives the desired claim.

\begin{rem}
The estimate \eqref{eq:packet-max} is the main reason
for our choice of $L^1$-normalized wave packets.
Heuristically, we may regard
$|\langle f,\varphi_{s,i}\rangle|$ as an average of $f$
over $I_s$.
\end{rem}

\bigskip 

\section{Reductions and the wave-packet model}\label{sec:reduction}
To prove Theorem~\ref{thm:main}, we need to show that
\begin{equation}\label{eq:level-target}
 \lambda\,|\{|\textnormal{\bf BHT}(f_1,f_2)|>\lambda\}|^{3/2}
 \lesssim_{p_1,p_2}\norm{f_1}_{p_1}\norm{f_2}_{p_2},
 \qquad \lambda>0.
\end{equation}
We begin with the following proposition, in which the second input
is bounded by the characteristic function of a set of finite measure.

\begin{prop}\label{prop:restricted}
Let $p_1,p_2$ satisfy \eqref{eq:critical}. For every
$f_1\in L^{p_1}(\R)$, every measurable set $F_2$ with
$0<|F_2|<\infty$, and every measurable function $f_2$ with
$|f_2|\le\one_{F_2}$, we have
\begin{equation}\label{eq:restricted-operator}
 \norm{\textnormal{\bf BHT}(f_1,f_2)}_{L^{2/3,\infty}}
 \lesssim_{p_1,p_2}\norm{f_1}_{p_1}|F_2|^{1/p_2}.
\end{equation}
\end{prop}

Sections~\ref{sec:reduction}--\ref{sec:model-proof} are devoted
to the proof of Proposition~\ref{prop:restricted}. The first step
is to reduce \eqref{eq:restricted-operator} to a trilinear estimate.
After normalizing the inputs and testing set, we pass to a
wave-packet model. Finally, 
only in Section~\ref{sec:completion}, we remove the restriction
$|f_2|\le\one_{F_2}$.

We first observe that, under the assumptions of
Proposition~\ref{prop:restricted}, $\textnormal{\bf BHT}(f_1,f_2)$
is well defined as an $L^1(\R)$ function. This follows from the classical ${\bf BHT}$ estimates of Lacey and Thiele~\cite{LT97,LT99}.

\subsection{Reduction to a trilinear estimate}\label{subsec:testing}
Let
\begin{equation}\label{eq:continuous-form}
 \Lambda(f_1,f_2,f_3)
 =\int_\R\textnormal{\bf BHT}(f_1,f_2)(x)f_3(x)\,dx.
\end{equation}
be the trilinear form associated to $\textnormal{\bf BHT}$. For a measurable set $F_3$ of finite positive measure, we call
$F_3'\subset F_3$ a \emph{major subset} if
$|F_3'|\ge |F_3|/2$. The following reduction is standard. 

\begin{lem}\label{lem:testing}
Fix $p_1,p_2$ satisfying \eqref{eq:critical}. Suppose that for every
$f_1\in L^{p_1}(\R)$ and measurable sets $F_2,F_3$ of finite
positive measure, we can choose a major subset $F_3'\subset F_3$
such that
\begin{equation}\label{eq:restricted-form}
 |\Lambda(f_1,f_2,f_3)|
 \lesssim_{p_1,p_2}\norm{f_1}_{p_1}
        |F_2|^{1/p_2}|F_3|^{-1/2}
\end{equation}
for all $|f_2|\le\one_{F_2}$ and $|f_3|\le\one_{F_3'}$.
Then Proposition~\ref{prop:restricted} holds.
\end{lem}

\begin{proof}
Fix $f_1,f_2$ as in Proposition~\ref{prop:restricted} and
$\lambda>0$. Let
\[
 F_3=\{x\in\R:|\textnormal{\bf BHT}(f_1,f_2)(x)|>\lambda\}.
\]
Since $\textnormal{\bf BHT}(f_1,f_2)\in L^1(\R)$, by Chebyshev, we have $
|F_3|\le \lambda^{-1}\norm{\textnormal{\bf BHT}(f_1,f_2)}_1<\infty$. 
If $|F_3|=0$, the desired level-set estimate is immediate.
Otherwise, choose $F_3'$ as in the hypothesis. Then taking
\[
 f_3=
 \frac{\overline{\textnormal{\bf BHT}(f_1,f_2)}}
 {|\textnormal{\bf BHT}(f_1,f_2)|}\one_{F_3'},
\]
with $f_3=0$ outside $F_3'$, we see that
\begin{align*}
 \frac{\lambda}{2}|F_3|
 &\le \lambda|F_3'|
 \le \int_{F_3'}|\textnormal{\bf BHT}(f_1,f_2)(x)|\,dx\\
 &=\Lambda(f_1,f_2,f_3)
 \lesssim_{p_1,p_2}\norm{f_1}_{p_1}|F_2|^{1/p_2}|F_3|^{-1/2}.
\end{align*}
Multiplying $|F_3|^{1/2}$ on both sides of the above estimate, we deduce that
\[
 \lambda\,|\{|\textnormal{\bf BHT}(f_1,f_2)|>\lambda\}|^{3/2}
 \lesssim_{p_1,p_2}\norm{f_1}_{p_1}|F_2|^{1/p_2}.
\]
This completes the proof. 
\end{proof}

\subsection{Normalization and the major subset}\label{subsec:normalization}

We now turn to the form estimate
\eqref{eq:restricted-form}. In the second reduction, we show that it suffices to prove the following normalized version.

\begin{prop}\label{prop:normalized}
Let $p_1,p_2$ satisfy \eqref{eq:critical}.
Suppose that $f_1\in L^{p_1}(\R)$ and that $F_2,F_3\subset\R$
are measurable sets satisfying
\begin{equation}\label{eq:normalized-inputs}
 \norm{f_1}_{p_1}=1,\qquad
 0<|F_2|<1,\qquad |F_3|=1.
\end{equation}
Then there exists a measurable subset $F_3'\subset F_3$ with
$|F_3'|\ge1/2$, depending only on $f_1,F_2,F_3$, such that
\begin{equation}\label{eq:normalized-form}
 |\Lambda(f_1,f_2,f_3)|
 \lesssim_{p_1,p_2}|F_2|^{1/p_2}
\end{equation}
for every $|f_2|\le\one_{F_2}$ and $|f_3|\le\one_{F_3'}$.
\end{prop}

\begin{proof}[Proof of \eqref{eq:restricted-form} assuming
Proposition~\ref{prop:normalized}]
Fix $f_1\in L^{p_1}(\R)$ and measurable sets $F_2,F_3$
of finite positive measure.
If $\norm{f_1}_{p_1}=0$, take $F_3'=F_3$.
We therefore assume that $\norm{f_1}_{p_1}>0$.

\vspace{0.1cm}

First suppose that $|F_2|\ge|F_3|$.
Take $F_3'=F_3$. By H\"older and the strong-type estimates
of {\bf BHT}, we have 
\begin{align*}
 |\Lambda(f_1,f_2,f_3)|
 &\lesssim_{p_1}
 \norm{f_1}_{p_1}\norm{f_2}_{2p_1'}\norm{f_3}_{2p_1'}\\
 &\lesssim_{p_1} \norm{f_1}_{p_1}
       |F_2|^{1/(2p_1')}|F_3|^{1/(2p_1')}\\
 &\lesssim_{p_1} \norm{f_1}_{p_1}
       |F_2|^{1/p_2}|F_3|^{-1/2}
\end{align*}
for all $|f_2|\le\one_{F_2}$ and $|f_3|\le\one_{F_3}$. Here, in the last estimate above, we used the facts that  $|F_3|\le|F_2|$ and
$1/p_2=1/p_1'+1/2$.

\vspace{0.1cm}

It remains to consider $0<|F_2|<|F_3|$.
Put $m=|F_3|$ and define
\[
 \widetilde f_1(x)
 :=\frac{m^{1/p_1}}{\norm{f_1}_{p_1}}f_1(mx),
\]
and
\[
 \widetilde F_i:=\{x\in\R:mx\in F_i\},
 \qquad i=2,3.
\]
Then
\[
 \norm{\widetilde f_1}_{p_1}=1,\qquad
 0<|\widetilde F_2|=\frac{|F_2|}{m}<1,\qquad
 |\widetilde F_3|=1.
\]
Applying now Proposition~\ref{prop:normalized} to
$\widetilde f_1,\widetilde F_2,\widetilde F_3$, we choose
$\widetilde F_3'\subset\widetilde F_3$ with
$|\widetilde F_3'|\ge1/2$. Note that 
the subset $\widetilde F_3'$ depends only on
$\widetilde f_1,\widetilde F_2,\widetilde F_3$.

Next, define
\[
 F_3':=\{mx:x\in\widetilde F_3'\}.
\]
It is clear that $F_3'$ depends only on $f_1,F_2,F_3$. Moreover,  $F_3'\subset F_3$ and $|F_3'|=m|\widetilde F_3'|\ge|F_3|/2$.

Finally, for any $f_2,f_3$ satisfying
$|f_2|\le\one_{F_2}$ and $|f_3|\le\one_{F_3'}$,
let
\[
 \widetilde f_i(x):=f_i(mx),\qquad i=2,3.
\]
Then
$|\widetilde f_2|\le\one_{\widetilde F_2}$ and
$|\widetilde f_3|\le\one_{\widetilde F_3'}$.
By the choice of $\widetilde F_3'$, we have
\[
 |\Lambda(\widetilde f_1,\widetilde f_2,\widetilde f_3)|
 \lesssim_{p_1,p_2}|\widetilde F_2|^{1/p_2}.
\]
By a change of variables, we deduce that
\begin{align*}
 |\Lambda(f_1,f_2,f_3)|
 &=\norm{f_1}_{p_1}m^{1-1/p_1}
   |\Lambda(\widetilde f_1,\widetilde f_2,\widetilde f_3)|\\
 &\lesssim_{p_1,p_2}
   \norm{f_1}_{p_1}m^{1-1/p_1}
   \left(\frac{|F_2|}{m}\right)^{1/p_2}\\
 &=\norm{f_1}_{p_1}|F_2|^{1/p_2}|F_3|^{-1/2},
\end{align*}
where the last equality follows from
$1/p_1+1/p_2=3/2$.

This proves \eqref{eq:restricted-form}.
\end{proof}

We now turn to the proof of Proposition~\ref{prop:normalized}.
Under the assumptions \eqref{eq:normalized-inputs}, we choose
$F_3'$ by removing an exceptional set from $F_3$.
For a sufficiently large constant $C$, set
\begin{equation}\label{eq:exceptional}
 \begin{split}
 E&=\{M_{p_1}f_1>C\}\cup\{M\one_{F_2}>C^{p_2}|F_2|\},\\
 \widetilde E&=\{M\one_E>1/10\},\qquad \textrm{and} \qquad 
 F_3'=F_3\setminus\widetilde E.
 \end{split}
\end{equation}
We record some properties for these sets. 

\begin{lem}\label{lem:major-subset}
Under the assumptions \eqref{eq:normalized-inputs}, the sets in
\eqref{eq:exceptional} satisfy
\[
 |E|\lesssim C^{-p_1}+C^{-p_2} \qquad \textrm{and} \qquad 
 |\widetilde E|\lesssim |E|.
\]
For $C$ sufficiently large, $|F_3'|\ge1/2$.
Moreover, if an interval $Q$ is contained in $E$, then
$3Q\subset\widetilde E$. In particular, every
$|f_3|\le\one_{F_3'}$ vanishes on $3Q$.
\end{lem}

\begin{proof}
For the first assertion, we apply the weak $(1,1)$ bound
for $M$ to obtain
\begin{align*}
 |E|
 &\le |\{M(|f_1|^{p_1})>C^{p_1}\}|
      +|\{M\one_{F_2}>C^{p_2}|F_2|\}|\\
 &\lesssim C^{-p_1}\norm{f_1}_{p_1}^{p_1}
      +(C^{p_2}|F_2|)^{-1}|F_2|\\
 &=C^{-p_1}+C^{-p_2}.
\end{align*}
Applying the same bound to $\one_E$, we also have
$|\widetilde E|
 =|\{M\one_E>1/10\}|
 \lesssim |E|$.
Therefore, for $C$ sufficiently large, $|\widetilde E|\le1/2$.
Since $|F_3|=1$, it follows that
\[
 |F_3'|=|F_3\setminus\widetilde E|
 \ge |F_3|-|\widetilde E|\ge\frac12.
\]

For the second assertion, let $Q\subset E$ be an interval.
For every $x\in3Q$, we have
\[
 M\one_E(x)
 \ge\frac1{|3Q|}\int_{3Q}\one_E(y)\,dy
 \ge\frac{|Q|}{|3Q|}
 >\frac1{10}.
\]
Thus $3Q\subset\widetilde E$.
Since $F_3'\subset\widetilde E^c$, every
$|f_3|\le\one_{F_3'}$ vanishes on $3Q$.
\end{proof}

The sets in \eqref{eq:exceptional} will be fixed throughout the rest of the proof of Proposition~\ref{prop:normalized}.
We remark that the set $F_3'$ depends only on $f_1,F_2,F_3$ and is independent
of the functions $f_2,f_3$.
It remains to prove \eqref{eq:normalized-form} for every
$|f_2|\le\one_{F_2}$ and $|f_3|\le\one_{F_3'}$.

\subsection{Reduction to a wave-packet model}\label{subsec:model}

We begin by applying the wave-packet discretization
in~\cite[Section~4]{DPT} to reduce
Proposition~\ref{prop:normalized} to uniform estimates
for finite wave-packet models.
To describe such models, we use the tiles and tri-tiles of
Definition~\ref{def:tile} and recall the notion of a
well-discretized collection of tri-tiles, following the presentation in~\cite[Section~4.2]{DPT}.

\begin{defn} \label{defn: well-discretized}
Fix a sufficiently large constant $R>10$.
A finite collection ${\bf S}$ of tri-tiles is called
\emph{well-discretized} if the following conditions hold:
\begin{enumerate}
\item[(i)]
The collections
\[
 \{I_s:s\in{\bf S}\},
 \qquad
 \{10\omega_{s,i}:s\in{\bf S}\},
 \quad i=1,2,3,
\]
are grids. Here, a collection of intervals is a \emph{grid} if any
two of its elements are disjoint or one contains the other.

\item[(ii)]
If $s,s'\in{\bf S}$ and $|I_s|<|I_{s'}|$, then
\[
 |I_s|\le R^{-10}|I_{s'}|.
\]

\item[(iii)]
If $s,s'\in{\bf S}$ are distinct and $I_s=I_{s'}$, then
\[
 \omega_{s,i}\cap\omega_{s',i}=\varnothing,
 \qquad i=1,2,3.
\]

\item[(iv)]
Suppose that $s,s'\in{\bf S}$ satisfy
$|I_{s'}|\le |I_s|$ and
\[
 2\omega_{s,j}\cap2\omega_{s',j}\ne\varnothing
\]
for some $j\in\{1,2,3\}$. Then
\[
 R\omega_{s,k}\subseteq R\omega_{s',k},
 \qquad k=1,2,3.
\]
If, in addition, $|I_{s'}|<|I_s|$, then
\[
 10\omega_{s,k}\cap10\omega_{s',k}=\varnothing,
 \qquad k\ne j.
\]
\end{enumerate}
\end{defn}

The constant $R$ is fixed throughout the proof and does \emph{not}
depend on the input functions.
 For each $s\in{\bf S}$ and $i=1,2,3$, let
$\varphi_{s,i}$ denote the $L^1$-normalized wave packet
associated with the tile $s_i$, as introduced in
Section~\ref{sec:prelim}. We require that the bounds in
\eqref{eq:packet}--\eqref{eq:packet-deriv} hold with constants
independent of $s$ and $i$.
For coefficients $\eps_s\in\C$ with $|\eps_s|\le1$, define
\begin{equation}\label{eq:model}
 \Lambda_{\bf S}(f_1,f_2,f_3)
 =\sum_{s\in{\bf S}}\eps_s|I_s|
       \prod_{i=1}^3\langle f_i,\varphi_{s,i}\rangle.
\end{equation}

\begin{rem} \label{rem:normalization}
The model form in~\cite[(4.3)]{DPT} is written using
$L^2$-normalized packets and therefore has the coefficient
$|I_s|^{-1/2}$ rather than $|I_s|$ in \eqref{eq:model}. To compare the two normalizations, write
\[
 \psi_{s,i}:=|I_s|^{1/2}\varphi_{s,i},
 \qquad i=1,2,3.
\]
Then the packets $\psi_{s,i}$ are the $L^2$-normalized, and hence
\[
 |I_s|^{-1/2}\prod_{i=1}^3
       \langle f_i,\psi_{s,i}\rangle
 =
 |I_s|\prod_{i=1}^3
       \langle f_i,\varphi_{s,i}\rangle.
\]
Thus \eqref{eq:model} coincides with the model form in \cite[(4.3)]{DPT}. 
\end{rem}

We now state the wave-packet reduction that allows us to estimate
$\Lambda$ using the finite model forms $\Lambda_{\bf S}$.

\begin{lem}\label{lem:model-reduction}
Assume \eqref{eq:normalized-inputs} and choose $F_3'$ as in
\eqref{eq:exceptional}. For every $|f_2|\le\one_{F_2}$ and
$|f_3|\le\one_{F_3'}$, we have
\begin{equation}\label{eq:discrete-reduction}
 |\Lambda(f_1,f_2,f_3)|
 \lesssim
 \sup_{\bf S}|\Lambda_{\bf S}(f_1,f_2,f_3)|
 +\left|\int_{\R}f_1(x)f_2(x)f_3(x)\,dx\right|.
\end{equation}
Here, the supremum is taken over all finite well-discretized
collections ${\bf S}$, all associated families of wave packets
satisfying \eqref{eq:packet}--\eqref{eq:packet-deriv}, and all
coefficients $|\eps_s|\le1$. 
The implicit constant in \eqref{eq:discrete-reduction} is
independent of the input functions and the number of tri-tiles.
\end{lem}

\begin{proof}
The lemma follows from~\cite[Section~4.3]{DPT}.
More precisely, one can argue exactly as in~\cite[(4.4)--(4.5)]{DPT},
with the $L^2$-normalized wave packets rescaled to the
$L^1$ normalization used here; see Remark \ref{rem:normalization}.
We leave the details to the interested reader.
\end{proof}

We next state the key model estimate that will imply
Proposition~\ref{prop:normalized}.

\begin{prop}\label{prop:model}
Let $p_1,p_2$ satisfy \eqref{eq:critical}, assume
\eqref{eq:normalized-inputs}, and choose $F_3'$ as in
\eqref{eq:exceptional}.
For every finite well-discretized collection ${\bf S}$,
we have
\begin{equation}\label{eq:model-target}
 |\Lambda_{\bf S}(f_1,f_2,f_3)|
 \lesssim_{p_1,p_2}|F_2|^{1/p_2}
\end{equation}
for all $|f_2|\le\one_{F_2}$ and $|f_3|\le\one_{F_3'}$.
The estimate is uniform over the collections ${\bf S}$,
the wave packets, and the coefficients specified in
Lemma~\ref{lem:model-reduction}.
\end{prop}

\begin{proof}[Proof of Proposition~\ref{prop:normalized}
assuming Proposition~\ref{prop:model}]
The set $F_3'$ defined in \eqref{eq:exceptional} depends only
on $f_1,F_2,F_3$, and Lemma~\ref{lem:major-subset} ensures that
$|F_3'|\ge1/2$.
By Lemma~\ref{lem:model-reduction} and
\eqref{eq:model-target}, it suffices to bound the product
integral in \eqref{eq:discrete-reduction}.

By the Lebesgue differentiation theorem, we have $E\subset\widetilde E$ up to a null set, and hence $F_3'\subset E^c$ up to a null set. Therefore, by the definition of $E$, we have 
\[
 |f_1(x)|\le M_{p_1}f_1(x)\le C
 \qquad\text{for almost every }x\in F_3'.
\]
Consequently, for every $|f_2|\le\one_{F_2}$ and
$|f_3|\le\one_{F_3'}$,
\begin{align*}
\left|\int_{\R}f_1(x)f_2(x)f_3(x)\,dx\right| \le \int_{F_2\cap F_3'}|f_1(x)|\,dx \le C|F_2|
 \le C|F_2|^{1/p_2},
\end{align*}
where the last inequality follows from $0<|F_2|<1$ and
$p_2>1$.
Combining this bound with \eqref{eq:discrete-reduction}
and \eqref{eq:model-target} gives
\eqref{eq:normalized-form}, which proves
Proposition~\ref{prop:normalized}.
\end{proof}

We now turn to the proof of Proposition~\ref{prop:model}.
We separate the tri-tiles in ${\bf S}$ into two collections: 
\begin{enumerate}
\item $I_s\subset E$;
\item  $I_s\not\subset E$.
\end{enumerate} 
For the first collection, Lemma~\ref{lem:major-subset}
ensures that $f_3$ vanishes on $3I_s$, allowing us to use
the decay of $\varphi_{s,3}$.
For the second collection, we organize the tri-tiles into
trees and prove a forest estimate.
We introduce the required tree structure and coefficient bounds in the next subsection. 

\subsection{Trees, sizes, and selection}\label{subsec:trees}

In this subsection, we introduce the tree structure, sizes,
and selection lemma needed for the forest estimate in
Proposition~\ref{prop:forest-p}.
The forest estimate will control the contribution of the
tri-tiles with $I_s\not\subset E$ in the proof of
Proposition~\ref{prop:model}.
We recall the definition of a tree, following the presentation in~\cite[Section~5.1]{DPT}.

\begin{defn}[Tree]\label{def:tree}
Let $j\in\{1,2,3\}$.
A \emph{$j$-tree} is a nonempty subcollection
${\bf T}\subset{\bf S}$ together with an interval
$I_{\bf T}$ and a frequency $\xi_{{\bf T},j}$ such that,
for every $s\in{\bf T}$,
\[
 I_s\subset I_{\bf T},\qquad
 \xi_{{\bf T},j}\in2\omega_{s,j}.
\]
We call $I_{\bf T}$ the \emph{top interval}.
A collection is called a \emph{tree} if it is a $j$-tree
for some $j\in\{1,2,3\}$.

For $i\in\{1,2,3\}$, the \emph{$i$th coordinate} of
${\bf T}$ is the collection obtained by taking the $i$th
tile from each tri-tile in ${\bf T}$, namely
\[
 \{s_i:s\in{\bf T}\}
 =
 \{I_s\times\omega_{s,i}:s\in{\bf T}\}.
\]
We say that ${\bf T}$ is \emph{lacunary in coordinate $i$}
if there exists a frequency $\xi\in\R$ such that,
for every $s\in{\bf T}$,
\[
 \xi\in R\omega_{s,i}\setminus2\omega_{s,i}.
\]
The frequency $\xi$ may depend on ${\bf T}$ and $i$,
but must be the same for every $s\in{\bf T}$.
\end{defn}

We next establish some properties of trees that will be
used in Section~\ref{sec:local}.

\begin{lem}\label{lem:tree-frequencies}
Let ${\bf T}$ be a $j$-tree with associated frequency
$\xi_{{\bf T},j}$.
Then, for each $i\ne j$, there exists a frequency
$\xi_{{\bf T},i}$ such that, for every $s\in{\bf T}$,
\[
 \xi_{{\bf T},i}\in
 R\omega_{s,i}\setminus2\omega_{s,i}.
\]
In particular, ${\bf T}$ is lacunary in every coordinate
$i\ne j$.

Moreover, with $K=(R+1)/2$, we have, for every
$s\in{\bf T}$ and $i\in\{1,2,3\}$,
\begin{equation}\label{eq:topfreq}
 \omega_{s,i}\subset
 \left[
   \xi_{{\bf T},i}-\frac{K}{|I_s|},
   \xi_{{\bf T},i}+\frac{K}{|I_s|}
 \right].
\end{equation}
\end{lem}

\begin{proof}
Choose $s_*\in{\bf T}$ with maximal spatial length.
For every $s\in{\bf T}$, we have
\[
 \xi_{{\bf T},j}
 \in2\omega_{s_*,j}\cap2\omega_{s,j}.
\]
By the well-discretized conditions, for every
$s\in{\bf T}$ and $i\in\{1,2,3\}$,
\[
 R\omega_{s_*,i}\subseteq R\omega_{s,i}.
\]
If $|I_s|=|I_{s_*}|$, then
$|R\omega_{s_*,i}|=|R\omega_{s,i}|$, and hence
\begin{equation}\label{eq:tree-equal-scale}
 \omega_{s,i}=\omega_{s_*,i}.
\end{equation}
If $|I_s|<|I_{s_*}|$, then for each $i\ne j$,
\begin{equation}\label{eq:tree-frequency-separation}
 10\omega_{s_*,i}\cap10\omega_{s,i}=\varnothing.
\end{equation}

Fix $i\ne j$ and choose
\[
 \xi_{{\bf T},i}
 \in3\omega_{s_*,i}\setminus2\omega_{s_*,i}.
\]
By \eqref{eq:tree-equal-scale} and
\eqref{eq:tree-frequency-separation}, we have
$\xi_{{\bf T},i}\notin2\omega_{s,i}$ for every $s\in{\bf T}$.
Moreover, since $R>10$, for every $s\in{\bf T}$, we have
\[
 \xi_{{\bf T},i}
 \in3\omega_{s_*,i}
 \subset R\omega_{s_*,i}
 \subseteq R\omega_{s,i}.
\]
Therefore, for every $s\in{\bf T}$, we see that
\[
 \xi_{{\bf T},i}
 \in R\omega_{s,i}\setminus2\omega_{s,i}.
\]
Thus the $i$th coordinate is lacunary. This proves the first assertion. 

\vspace{0.1cm}

Next, we show the second assertion. For every $s\in{\bf T}$ and $i\in\{1,2,3\}$, the chosen
frequency $\xi_{{\bf T},i}$ belongs to $R\omega_{s,i}$.
Therefore, for every $\eta\in\omega_{s,i}$,
\[
 |\eta-\xi_{{\bf T},i}|
 \le\frac{R+1}{2}|\omega_{s,i}|
 =\frac{K}{|I_s|}.
\]
This proves \eqref{eq:topfreq}.
\end{proof}

\begin{lem}\label{lem:lacunary-disjointness}
Suppose that the $i$th coordinate of a tree ${\bf T}$
is lacunary.
For every $s,s'\in{\bf T}$, we have
\begin{equation}\label{eq:lacunary}
 \omega_{s,i}\ne\omega_{s',i}
 \quad\Longrightarrow\quad
 \omega_{s,i}\cap\omega_{s',i}=\varnothing.
\end{equation}
\end{lem}

\begin{proof}
Choose a frequency $\xi$ such that, for every $s\in{\bf T}$,
\[
 \xi\in R\omega_{s,i}\setminus2\omega_{s,i}.
\]
Fix $s,s'\in{\bf T}$ and suppose that
$\omega_{s,i}\cap\omega_{s',i}\ne\varnothing$.

If $|\omega_{s,i}|=|\omega_{s',i}|$, then $10\omega_{s,i}\cap10\omega_{s',i}\ne\varnothing$ and $|10\omega_{s,i}|=|10\omega_{s',i}|$. By the well-discretized condition, we have $10\omega_{s,i}=10\omega_{s',i}$, and hence $\omega_{s,i}=\omega_{s',i}$.

If the lengths differ, we may assume
$|\omega_{s,i}|<|\omega_{s',i}|$.
Using the well-discretized condition again,
\[
 |\omega_{s,i}|\le R^{-10}|\omega_{s',i}|.
\]
Since $\omega_{s,i}$ and $\omega_{s',i}$ intersect, we have
\begin{align*}
 |\xi-c(\omega_{s',i})|
 &\le
 |\xi-c(\omega_{s,i})|
 +|c(\omega_{s,i})-c(\omega_{s',i})|\\
 &\le
 \frac{R}{2}|\omega_{s,i}|
 +\frac{|\omega_{s,i}|+|\omega_{s',i}|}{2}\\
 &\le
 \left(\frac12+\frac{R+1}{2R^{10}}\right)|\omega_{s',i}|
 <|\omega_{s',i}|.
\end{align*}
Hence $\xi\in2\omega_{s',i}$, however, this contradicts with the choice of $\xi$.
This proves \eqref{eq:lacunary}.
\end{proof}

\begin{lem}\label{lem:tree-intervals}
Let ${\bf T}$ be a tree.
For every $\ell>0$, distinct spatial intervals of length
$\ell$ occurring in ${\bf T}$ are disjoint, and their
centers are separated by at least $\ell$.
Moreover, if $\ell<\ell'$ are two spatial lengths occurring
in ${\bf T}$, then
\[
 \ell\le R^{-10}\ell'.
\]
\end{lem}

\begin{proof}
For the first assertion, note that the spatial intervals
form a grid, so any two intersecting intervals of equal
length must coincide.
Hence distinct spatial intervals of length $\ell$ are
disjoint, and their centers are separated by at least $\ell$.
The second assertion follows directly from
Definition~\ref{defn: well-discretized}.
\end{proof}

\begin{lem}\label{lem:tree-counting}
Let ${\bf T}$ be a tree.
For every interval $I$, we have
\[
 \#\{s\in{\bf T}:I_s=I\}\lesssim1.
\]
Consequently, for every $\ell>0$,
\[
 \#\{s\in{\bf T}:|I_s|=\ell\}
 \lesssim\frac{|I_{\bf T}|}{\ell}.
\]
\end{lem}

\begin{proof}
For the first assertion, suppose that ${\bf T}$ is a
$j$-tree, and fix an interval $I$ with $|I|=\ell$.
The frequency intervals $\omega_{s,j}$ corresponding to
distinct tri-tiles with $I_s=I$ are disjoint and have
length $\ell^{-1}$.
For each such tri-tile, the condition
$\xi_{{\bf T},j}\in2\omega_{s,j}$ implies that
\[
 \omega_{s,j}\subset
 \left[
   \xi_{{\bf T},j}-\frac{3}{2\ell},
   \xi_{{\bf T},j}+\frac{3}{2\ell}
 \right].
\]
Since the interval on the right has length $3\ell^{-1}$, it contains at most three disjoint intervals of length
$\ell^{-1}$.
This proves the first assertion.

For the second assertion, Lemma~\ref{lem:tree-intervals}
shows that distinct spatial intervals of length $\ell$
are disjoint.
Since all of these intervals are contained in $I_{\bf T}$,
their number is at most $|I_{\bf T}|/\ell$.
Each interval corresponds to at most three tri-tiles
by the first assertion.
The desired bound follows.
\end{proof}

We next introduce the notion of a forest.

\begin{defn}[Forest]\label{def:forest}
Let ${\bf S}$ be a finite well-discretized collection
of tri-tiles.
A \emph{forest} in ${\bf S}$ is a finite family
$\mathcal F$ of trees ${\bf T}\subset{\bf S}$ that are
pairwise disjoint as sets of tri-tiles.
We call $\mathcal F$ a \emph{forest decomposition}
of ${\bf S}$ if
\[
 {\bf S}=\bigsqcup_{{\bf T}\in\mathcal F}{\bf T}.
\]
\end{defn}

For a finite well-discretized collection ${\bf S}$, we define the size of a function associated to ${\bf S}$ in each coordinate.

\begin{defn}[Size]\label{def:size}
Let ${\bf S}$ be a finite well-discretized collection
of tri-tiles, and fix $i\in\{1,2,3\}$.
For $f\in L^1(\R)+L^\infty(\R)$, define
\begin{equation}\label{eq:size}
 \operatorname{size}_i(f;{\bf S})
 =
 \sup_{\substack{
   {\bf T}\subset{\bf S}\text{ is a tree}\\
   {\bf T}\text{ is lacunary in coordinate }i}}
 \ \sup_{\{\varphi_{s,i}\}_{s\in{\bf T}}}
 \left(
   \frac1{|I_{\bf T}|}
   \sum_{s\in{\bf T}}|I_s|\,
       |\langle f,\varphi_{s,i}\rangle|^2
 \right)^{1/2}.
\end{equation}
The inner supremum is taken over all families of
$L^1$-normalized wave packets $\varphi_{s,i}$ adapted to
$I_s\times\omega_{s,i}$ and satisfying
\eqref{eq:packet}--\eqref{eq:packet-deriv}.
Here, if ${\bf S}=\varnothing$, we define $\operatorname{size}_i(f;\varnothing)=0$.
\end{defn}

Here are some direct consequences of Definition~\ref{def:size}.
First, note that for every $s\in{\bf S}$ and $i\in\{1,2,3\}$, one has $|\langle f_i,\varphi_{s,i}\rangle|
 \le \operatorname{size}_i(f_i;{\bf S})$. 
This follows by taking ${\bf T}=\{s\}$ in \eqref{eq:size}.  More generally, for any ${\bf T}\subset{\bf S}$
that is lacunary in coordinate $i$, one has
\begin{equation}\label{eq:size-coefficients}
 \sum_{s\in{\bf T}}|I_s|\,
       |\langle f_i,\varphi_{s,i}\rangle|^2
 \le
 |I_{\bf T}|\operatorname{size}_i(f_i;{\bf S})^2.
\end{equation}
Note also that restricting ${\bf S}$ to a subcollection does not
increase its size, since the supremum is then taken over
a smaller collection of trees.

We now estimate the model form restricted to a single tree.

\begin{lem}[Tree estimate]\label{lem:tree}
Let ${\bf T}\subset{\bf S}$ be a tree.
Using the wave packets and coefficients in
\eqref{eq:model}, define
\[
 \Lambda_{\bf T}(f_1,f_2,f_3)
 :=\sum_{s\in{\bf T}}\eps_s|I_s|
 \prod_{i=1}^3\langle f_i,\varphi_{s,i}\rangle.
\]
Then
\begin{equation}\label{eq:tree}
 \begin{split}
 |\Lambda_{\bf T}(f_1,f_2,f_3)|
 &\le
 \sum_{s\in{\bf T}}|I_s|
 \prod_{i=1}^3|\langle f_i,\varphi_{s,i}\rangle|\\
 &\lesssim
 |I_{\bf T}|\prod_{i=1}^3
          \operatorname{size}_i(f_i;{\bf S}).
 \end{split}
\end{equation}
\end{lem}

\begin{proof}
Suppose that ${\bf T}$ is a $j$-tree.
By Lemma~\ref{lem:tree-frequencies}, ${\bf T}$ is lacunary
in both coordinates $i\ne j$.
Applying Cauchy--Schwarz in these coordinates, we have
\begin{align*}
\sum_{s\in{\bf T}}|I_s|
       \prod_{i=1}^3|\langle f_i,\varphi_{s,i}\rangle|
       &\le
 \sup_{s\in{\bf T}}|\langle f_j,\varphi_{s,j}\rangle|
 \prod_{i\ne j}
 \left(
   \sum_{s\in{\bf T}}|I_s|\,
       |\langle f_i,\varphi_{s,i}\rangle|^2
 \right)^{1/2}\\
 & \le
 |I_{\bf T}|\prod_{i=1}^3
          \operatorname{size}_i(f_i;{\bf S}).
\end{align*}
Here, in the last estimate above, we used 
\eqref{eq:size-coefficients}.
\end{proof}

To sum the estimates over a forest, we need to control
the sizes and the sum of the lengths of top intervals.
We record both estimates in the following lemma.

\begin{lem}\label{lem:size-selection}
Let ${\bf S}$ be a finite well-discretized collection
of tri-tiles as in Definition~\ref{defn: well-discretized}.
For every $f\in L^1(\R)+L^\infty(\R)$ and
$i\in\{1,2,3\}$, we have
\begin{equation}\label{eq:size-max}
 \operatorname{size}_i(f;{\bf S})
 \lesssim
 \sup_{s\in{\bf S}}\inf_{x\in I_s}Mf(x).
\end{equation}
If, in addition, $f\in L^2(\R)$ and
$\operatorname{size}_i(f;{\bf S})\le\sigma$ for some
$\sigma>0$, then there exist a subcollection
${\bf S}_{\mathrm{low}}$ and a forest $\mathcal F$ such that
\begin{equation}\label{eq:selection}
 \begin{split}
 {\bf S}
 &={\bf S}_{\mathrm{low}}
   \sqcup\bigsqcup_{{\bf T}\in\mathcal F}{\bf T},\\
 \operatorname{size}_i(f;{\bf S}_{\mathrm{low}})
 &\le\frac{\sigma}{2},
 \qquad \textrm{and} \qquad 
 \sum_{{\bf T}\in\mathcal F}|I_{\bf T}|
 \lesssim\sigma^{-2}\norm{f}_2^2.
 \end{split}
\end{equation}
The implicit constants depend only on $R$ and the fixed
wave-packet bounds.
\end{lem}

\begin{proof}
Fix $i\in\{1,2,3\}$.
For the first assertion, let ${\bf T}$ be a tree that is
lacunary in coordinate $i$, and let
$\{\varphi_{s,i}\}_{s\in{\bf T}}$ be any admissible family
in \eqref{eq:size}.
Set
\begin{equation} \label{eq: L2--normalizationV1}
 \psi_{s,i}:=|I_s|^{1/2}\varphi_{s,i}.
\end{equation}
Then
\[
 \frac1{|I_{\bf T}|}
 \sum_{s\in{\bf T}}|I_s|\,
       |\langle f,\varphi_{s,i}\rangle|^2
 =
 \frac1{|I_{\bf T}|}
 \sum_{s\in{\bf T}}|\langle f,\psi_{s,i}\rangle|^2.
\]

Choose a frequency $\xi$ satisfying the lacunarity condition
in Definition~\ref{def:tree}.
For every $s\in{\bf T}$, we then have
\begin{equation}\label{eq:packet-lacunarity}
 |\xi-c(\omega_{s,i})|\le\frac{R}{2|I_s|},
 \qquad
 \widehat{\psi}_{s,i}(\xi)=0.
\end{equation}
Using \eqref{eq:packet-lacunarity} and
\eqref{eq:packet}--\eqref{eq:packet-deriv}, we check directly
that the packets $\{\psi_{s,i}\}_{s\in{\bf T}}$ satisfy
the assumptions of~\cite[Lemma~3.5]{DPT}.
Applying that lemma to $\{\psi_{s,i}\}_{s\in{\bf T}}$,
as in the argument of~\cite[Section~5.1]{DPT},
we have
\[
 \left(
 \frac1{|I_{\bf T}|}
 \sum_{s\in{\bf T}}|\langle f,\psi_{s,i}\rangle|^2
 \right)^{1/2}
 \lesssim
 \sup_{s\in{\bf T}}\inf_{x\in I_s}Mf(x).
\]
Taking the suprema over the admissible packet families
and the trees, we conclude \eqref{eq:size-max}.

For the second assertion, we choose an auxiliary family
$\{\varphi_{s,i}\}_{s\in{\bf S}}$, depending on $f$, such
that, for every $s\in{\bf S}$,
\[
 \sup_{\varphi}|\langle f,\varphi\rangle|
 \le2|\langle f,\varphi_{s,i}\rangle|,
\]
where $\varphi$ ranges over all admissible packets associated
with $I_s\times\omega_{s,i}$.
We keep this family fixed on ${\bf S}$ throughout the
rest of the proof.

For every nonempty subcollection ${\bf S}'\subset{\bf S}$,
define
\[
 \operatorname{size}_i^\varphi(f;{\bf S}')
 =
 \sup_{\substack{
   {\bf T}\subset{\bf S}'\text{ is a tree}\\
   {\bf T}\text{ is lacunary in coordinate }i}}
 \left(
   \frac1{|I_{\bf T}|}
   \sum_{s\in{\bf T}}|I_s|\,
       |\langle f,\varphi_{s,i}\rangle|^2
 \right)^{1/2},
\]
where the packets $\varphi_{s,i}$ belong to the fixed
family chosen on ${\bf S}$.
We also set $\operatorname{size}_i^\varphi(f;\varnothing)=0$.
By the choice of the packets, for every subcollection
${\bf S}'\subset{\bf S}$, we have
\begin{equation}\label{eq:size-fixed-comparison}
 \operatorname{size}_i^\varphi(f;{\bf S}')
 \le\operatorname{size}_i(f;{\bf S}')
 \le2\operatorname{size}_i^\varphi(f;{\bf S}').
\end{equation}
Using the same normalization as in \eqref{eq: L2--normalizationV1}, 
we can write, for every nonempty subcollection
${\bf S}'\subset{\bf S}$,
\[
 \operatorname{size}_i^\varphi(f;{\bf S}')
 =
 \sup_{\substack{
   {\bf T}\subset{\bf S}'\text{ is a tree}\\
   {\bf T}\text{ is lacunary in coordinate }i}}
 \left(
   \frac1{|I_{\bf T}|}
   \sum_{s\in{\bf T}}|\langle f,\psi_{s,i}\rangle|^2
 \right)^{1/2}.
\]
By \eqref{eq:size-fixed-comparison}, we have
$\operatorname{size}_i^\varphi(f;{\bf S})\le\sigma$.
Arguing as in the proof of~\cite[Lemma~7.7]{MTT}
for the fixed family $\{\psi_{s,i}\}_{s\in{\bf S}}$,
we obtain a subcollection ${\bf S}^{(1)}\subset{\bf S}$
and a forest $\mathcal F_1$ such that
\[
 {\bf S}
 =
 {\bf S}^{(1)}
 \sqcup\bigsqcup_{{\bf T}\in\mathcal F_1}{\bf T},
\]
with
\[
 \operatorname{size}_i^\varphi(f;{\bf S}^{(1)})
 \le\frac{\sigma}{2},
 \qquad\text{and}\qquad
 \sum_{{\bf T}\in\mathcal F_1}|I_{\bf T}|
 \lesssim\sigma^{-2}\norm{f}_2^2.
\]

We now repeat the same argument on ${\bf S}^{(1)}$, with
$\sigma/2$ in place of $\sigma$, to deduce a subcollection
${\bf S}_{\mathrm{low}}$ and a forest $\mathcal F_2$ such that
\[
 {\bf S}^{(1)}
 =
 {\bf S}_{\mathrm{low}}
 \sqcup\bigsqcup_{{\bf T}\in\mathcal F_2}{\bf T},
\]
with
\[
 \operatorname{size}_i^\varphi(f;{\bf S}_{\mathrm{low}})
 \le\frac{\sigma}{4},
 \qquad \textrm{and} \qquad 
 \sum_{{\bf T}\in\mathcal F_2}|I_{\bf T}|
 \lesssim(\sigma/2)^{-2}\norm{f}_2^2.
\]

Set $\mathcal F:=\mathcal F_1\cup\mathcal F_2$.
The trees in $\mathcal F_2$ are contained in
${\bf S}^{(1)}$ and are therefore disjoint from those
in $\mathcal F_1$ as sets of tri-tiles.
Thus $\mathcal F$ is a forest, and
\[
 {\bf S}
 =
 {\bf S}_{\mathrm{low}}
 \sqcup\bigsqcup_{{\bf T}\in\mathcal F}{\bf T}.
\]
Adding the bounds for the top interval lengths in
$\mathcal F_1$ and $\mathcal F_2$, we have
\[
 \sum_{{\bf T}\in\mathcal F}|I_{\bf T}|
 \lesssim
 \left(\sigma^{-2}+(\sigma/2)^{-2}\right)\norm{f}_2^2
 \lesssim\sigma^{-2}\norm{f}_2^2.
\]
Finally, applying \eqref{eq:size-fixed-comparison}
to ${\bf S}_{\mathrm{low}}$, we deduce that
\[
 \operatorname{size}_i(f;{\bf S}_{\mathrm{low}})
 \le
 2\operatorname{size}_i^\varphi(f;{\bf S}_{\mathrm{low}})
 \le\frac{\sigma}{2}.
\]
This proves \eqref{eq:selection}.
\end{proof}

\subsection{The forest estimate needed for the model}

Under the assumptions of Proposition~\ref{prop:model},
let ${\bf S}$ be a finite well-discretized collection
such that $I_s\not\subset E$ for every $s\in{\bf S}$.
By \eqref{eq:exceptional} and \eqref{eq:size-max}, we have
\[
 \operatorname{size}_2(f_2;{\bf S})\lesssim |F_2|,
 \qquad
 \operatorname{size}_3(f_3;{\bf S})\lesssim1.
\]
We will use \eqref{eq:selection} to control the sum of the lengths of
top intervals at each size level.
To retain the $L^{p_1}$ information on the first input,
we need the following forest estimate, which is central to our argument. 

\begin{prop}[Forest estimate]\label{prop:forest-p}
Let $1<p<2$ and let $f\in L^p(\R)$ satisfy
$\norm{f}_p\le1$.
Fix $C_0>1$ and set $E_f=\{M_pf>C_0\}$.
Suppose that ${\bf S}$ is a finite well-discretized
collection of tri-tiles satisfying $I_s\not\subset E_f$
for every $s\in{\bf S}$, and let
\[
 {\bf S}=\bigsqcup_{{\bf T}\in\mathcal F}{\bf T}
\]
be a forest decomposition.
Write $N=\sum_{{\bf T}\in\mathcal F}|I_{\bf T}|$.
If $f_2\in L^1(\R)+L^\infty(\R)$ and $f_3\in L^2(\R)$
satisfy
$\operatorname{size}_2(f_2;{\bf S})<\infty$ and
$\operatorname{size}_3(f_3;{\bf S})<\infty$,
then
\begin{equation}\label{eq:forest-p}
 |\Lambda_{\bf S}(f,f_2,f_3)|
 \lesssim_{p,C_0}
 (N^{1/2}+N^{1/p})
 \operatorname{size}_2(f_2;{\bf S})
 \operatorname{size}_3(f_3;{\bf S}).
\end{equation}
\end{prop}

\begin{rem}
Later in the proof of Proposition~\ref{prop:model}, we will apply
Proposition~\ref{prop:forest-p} with
$f=f_1$, $p=p_1$, and $C_0=C$,
where $p_1$ is as in \eqref{eq:critical} and $C$ is
the fixed constant chosen in \eqref{eq:exceptional}.
\end{rem}

\bigskip

\section{The local estimate for a forest}\label{sec:local}

Our goal in this section is to prove Proposition~\ref{prop:forest-p}.

\subsection{The packet sum and its \texorpdfstring{$L^2$}{L2} bound}

Let ${\bf S}$ be a finite well-discretized collection
with a forest decomposition
\[
 {\bf S}=\bigsqcup_{{\bf T}\in\mathcal F}{\bf T},
 \qquad \textrm{and} \qquad 
 N=\sum_{{\bf T}\in\mathcal F}|I_{\bf T}|.
\]
Observe that the top intervals may overlap, and the sum defining $N$
includes one term for each tree.
For the moment, assume that
\[
 \operatorname{size}_2(f_2;{\bf S})\le1,
 \qquad
 \operatorname{size}_3(f_3;{\bf S})\le1.
\]
For the packets and coefficients in \eqref{eq:model}, define
\begin{equation}\label{eq:synthesis}
 \beta_s:=\overline{\eps_s|I_s|
       \langle f_2,\varphi_{s,2}\rangle
       \langle f_3,\varphi_{s,3}\rangle}, \qquad \textrm{and} \qquad 
 \Psi
 :=\sum_{s\in{\bf S}}\beta_s\varphi_{s,1}.
\end{equation}
Observe that 
\begin{equation}\label{eq:forest-pairing}
 \Lambda_{\bf S}(f,f_2,f_3)=\langle f,\Psi\rangle.
\end{equation}
For every tree ${\bf T}\subset{\bf S}$, we also define
\[
 \Psi_{\bf T}:=\sum_{s\in{\bf T}}\beta_s\varphi_{s,1}.
\]
We begin with an $L^2$ estimate for the pairing $\Lambda_{\bf S}(f, f_2, f_3)$.

\begin{lem}\label{lem:forest2}
Suppose that
$\operatorname{size}_2(f_2;{\bf S})\le1$ and
$\operatorname{size}_3(f_3;{\bf S})\le1$.
For every $h\in L^2(\R)$, we have
\begin{equation}\label{eq:forest-l2}
 |\Lambda_{\bf S}(h,f_2,f_3)|
 \lesssim N^{1/2}\norm{h}_2.
\end{equation}
Moreover, for every $s\in{\bf S}$, the coefficients and
packet sum satisfy
\begin{equation}\label{eq:coeff}
 |\beta_s|\le|I_s|,
 \qquad
 \sum_{s\in{\bf S}}\frac{|\beta_s|^2}{|I_s|}
 \lesssim N,
 \qquad
 \norm{\Psi}_2\lesssim N^{1/2}.
\end{equation}
\end{lem}

\begin{proof}
For the first assertion, the cases ${\bf S}=\varnothing$
and $h=0$ are immediate.
Otherwise, by homogeneity we may assume $\norm{h}_2=1$.
Set $\tau=N^{-1/2}$.
Since $h\in L^2(\R)$ and ${\bf S}$ is finite,
\eqref{eq:size-max} 
implies that $\operatorname{size}_1(h;{\bf S})<\infty$.
Hence, we can choose an integer $M\ge1$ such that
\[
 \operatorname{size}_1(h;{\bf S})\le2^M\tau.
\]
We apply \eqref{eq:selection} successively to the
remaining collection at thresholds
$2^M\tau,2^{M-1}\tau,\ldots,2\tau$.
This gives a subcollection ${\bf S}_{\mathrm{low}}$
and forests $\mathcal F_m$, indexed by $1\le m\le M$,
such that
\begin{equation}\label{eq:forest-selection-decomposition}
 {\bf S}
 =
 {\bf S}_{\mathrm{low}}
 \sqcup\bigsqcup_{m=1}^M
       \bigsqcup_{{\bf T}\in\mathcal F_m}{\bf T},
 \qquad \textrm{and} \qquad
 \operatorname{size}_1(h;{\bf S}_{\mathrm{low}})\le\tau.
\end{equation}
For every $1\le m\le M$ with $\mathcal F_m\ne\varnothing$,
we have
\begin{equation}\label{eq:forest-level-selection}
 \sup_{{\bf T}\in\mathcal F_m}
       \operatorname{size}_1(h;{\bf T})\le2^m\tau,
 \qquad \textrm{and} \qquad 
 \sum_{{\bf T}\in\mathcal F_m}|I_{\bf T}|
 \lesssim(2^m\tau)^{-2}.
\end{equation}
Moreover, we have
\begin{equation}\label{eq:forest-other-sizes}
 \sup_{{\bf T}\in\mathcal F_m}
 \operatorname{size}_2(f_2;{\bf T})
 \le \operatorname{size}_2(f_2;{\bf S})\le1,
 \qquad \textrm{and} \qquad 
 \sup_{{\bf T}\in\mathcal F_m}
 \operatorname{size}_3(f_3;{\bf T})
 \le \operatorname{size}_3(f_3;{\bf S})\le1,
\end{equation}
since every ${\bf T}\in\mathcal F_m$ is contained in
${\bf S}$ and size cannot increase under restriction.

As a consequence of the decomposition \eqref{eq:forest-selection-decomposition}, we may write
\begin{equation}\label{eq:forest-form-decomposition}
 \Lambda_{\bf S}(h,f_2,f_3)
 =
 \sum_{m=1}^M\sum_{{\bf T}\in\mathcal F_m}
       \Lambda_{\bf T}(h,f_2,f_3)
 +
 \Lambda_{{\bf S}_{\mathrm{low}}}(h,f_2,f_3).
\end{equation}
We now estimate the above two terms separately.

\smallskip
\noindent $\bullet$ \textsf{Estimate of the selected trees.}
For every $1\le m\le M$ with $\mathcal F_m\ne\varnothing$,
we apply \eqref{eq:tree}, together with
\eqref{eq:forest-level-selection} and
\eqref{eq:forest-other-sizes}, to obtain
\[
 \sum_{{\bf T}\in\mathcal F_m}
       |\Lambda_{\bf T}(h,f_2,f_3)|
 \lesssim
 2^m\tau\sum_{{\bf T}\in\mathcal F_m}|I_{\bf T}|
 \lesssim(2^m\tau)^{-1}.
\]
Summing over $1\le m\le M$, we have
\begin{equation}\label{eq:forest-selected}
 \sum_{m=1}^M\sum_{{\bf T}\in\mathcal F_m}
       |\Lambda_{\bf T}(h,f_2,f_3)|
 \lesssim
 \tau^{-1}\sum_{m=1}^M2^{-m}
 \le\tau^{-1}.
\end{equation}

\smallskip
\noindent $\bullet$ \textsf{Estimate of the remaining collection.}
Using the original forest $\mathcal F$, we decompose
the remaining collection as
\[
 {\bf S}_{\mathrm{low}}
 =
 \bigsqcup_{{\bf T}\in\mathcal F}
       ({\bf T}\cap{\bf S}_{\mathrm{low}}).
\]
Note that by Definition~\ref{def:tree}, we have if ${\bf T}\in\mathcal F$ is a $j$-tree and
${\bf T}\cap{\bf S}_{\mathrm{low}}\ne\varnothing$,
then ${\bf T}\cap{\bf S}_{\mathrm{low}}$ is also a
$j$-tree with top interval $I_{\bf T}$ and frequency
$\xi_{{\bf T},j}$.

Since ${\bf S}_{\mathrm{low}}\subset{\bf S}$, we have,
for every $i\in\{2,3\}$,
\[
 \operatorname{size}_i(f_i;{\bf S}_{\mathrm{low}})
 \le \operatorname{size}_i(f_i;{\bf S})\le1,
\]
and hence, by \eqref{eq:forest-selection-decomposition}
and \eqref{eq:tree}, we obtain
\[
 |\Lambda_{{\bf T}\cap{\bf S}_{\mathrm{low}}}(h,f_2,f_3)|
 \lesssim\tau|I_{\bf T}|.
\]
Summing over ${\bf T}\in\mathcal F$, we conclude that
\begin{equation}\label{eq:forest-low}
 |\Lambda_{{\bf S}_{\mathrm{low}}}(h,f_2,f_3)|
 \lesssim
 \tau\sum_{{\bf T}\in\mathcal F}|I_{\bf T}|
 =N\tau.
\end{equation}

Combining \eqref{eq:forest-selected} and
\eqref{eq:forest-low} with
\eqref{eq:forest-form-decomposition}, 
\[
 |\Lambda_{\bf S}(h,f_2,f_3)|
 \lesssim\tau^{-1}+N\tau=2N^{1/2},
\]
which proves the first assertion. 

\vspace{0.1cm}

Next, we prove the second assertion.
Applying Definition~\ref{def:size} to singleton trees,
we have, for every $s\in{\bf S}$ and $i\in\{2,3\}$,
\begin{equation}\label{eq:forest-singleton-size}
 |\langle f_i,\varphi_{s,i}\rangle|
 \le\operatorname{size}_i(f_i;{\bf S})\le1.
\end{equation}
Consequently, for every $s\in{\bf S}$,
\[
 |\beta_s|
 \le |I_s|\,
       |\langle f_2,\varphi_{s,2}\rangle|
       |\langle f_3,\varphi_{s,3}\rangle|
 \le |I_s|.
\]

Fix a $j$-tree ${\bf T}\in\mathcal F$ and choose
$i\in\{2,3\}\setminus\{j\}$.
By Lemma~\ref{lem:tree-frequencies}, ${\bf T}$ is
lacunary in coordinate $i$.
Applying \eqref{eq:size-coefficients} in coordinate $i$
and \eqref{eq:forest-singleton-size} in the other coordinate
in $\{2,3\}$, we obtain 
\[
 \sum_{s\in{\bf T}}\frac{|\beta_s|^2}{|I_s|}
 \le
 \sum_{s\in{\bf T}}|I_s|\,
       |\langle f_i,\varphi_{s,i}\rangle|^2
 \le |I_{\bf T}|.
\]
Summing over ${\bf T}\in\mathcal F$, 
\[
 \sum_{s\in{\bf S}}\frac{|\beta_s|^2}{|I_s|}
 \le\sum_{{\bf T}\in\mathcal F}|I_{\bf T}|=N,
\]
which proves the second bound in \eqref{eq:coeff}.

Finally, by \eqref{eq:forest-pairing},
\eqref{eq:forest-l2}, duality, we conclude that 
\[
 \norm{\Psi}_2
 =\sup_{\norm{h}_2=1}|\langle h,\Psi\rangle|
 \lesssim N^{1/2}.
\]
The proof of Lemma \ref{lem:forest2} is complete. 
\end{proof}

Note that in Proposition~\ref{prop:forest-p}, the input function $f$
belongs to $L^2$ on the complement of
$E_f=\{M_pf>C_0\}$.
We will therefore use \eqref{eq:forest-l2} to estimate
the pairing over $E_f^c$.
For the pairing over $E_f$, we need an additional bound
for $\Psi$ on the Whitney intervals of $E_f$.

\subsection{Local suprema on Whitney intervals}

We state the local supremum estimate for an arbitrary open set
$\Omega$ of finite measure.
In the proof of Proposition~\ref{prop:forest-p},
we will take $\Omega=E_f$.

\begin{prop}\label{prop:local}
Let $\Omega\subset\R$ be open with finite measure.
Fix a dyadic lattice $\calD$, and let $\mathcal W$
consist of the maximal intervals $H\in\calD$
such that $9H\subset\Omega$.
Suppose that ${\bf S}$ is a finite well-discretized
collection with a forest decomposition
\[
 {\bf S}=\bigsqcup_{{\bf T}\in\mathcal F}{\bf T},
\]
and that
$\operatorname{size}_2(f_2;{\bf S})\le1$ and
$\operatorname{size}_3(f_3;{\bf S})\le1$.
Assume that, for every $s\in{\bf S}$,
$I_s\not\subset\Omega$.
Then the packet sum in \eqref{eq:synthesis} satisfies
\begin{equation}\label{eq:local}
 \sum_{H\in\mathcal W}|H|\norm{\Psi}_{L^\infty(H)}
 \lesssim
 \sum_{{\bf T}\in\mathcal F}|I_{\bf T}|.
\end{equation}
\end{prop}

The local supremum estimate in Proposition~\ref{prop:local}
is the main ingredient in this paper.
We organize the proof into several lemmas.

\begin{lem}\label{lem:whitney}
The intervals in $\mathcal W$ are pairwise disjoint
and cover $\Omega$ up to endpoints.
For every $H\in\mathcal W$ and every tri-tile $s$
with $I_s\not\subset\Omega$, we have
\begin{equation}\label{eq:whitney-separation}
 |I_s|<|H|
 \quad\Longrightarrow\quad
 \dist(H,I_s)\ge3|H|.
\end{equation}
\end{lem}

\begin{proof}
The first assertion follows from the fact that the dyadic
intervals $H$ with $9H\subset\Omega$ cover $\Omega$ and
have bounded lengths, together with the nesting property
of the dyadic lattice.

For the second assertion, fix $H\in\mathcal W$ and
choose $z\in I_s\setminus\Omega$.
Since $9H\subset\Omega$, the point $z$ lies outside
$9H$, so $\dist(z,H)\ge4|H|$.
Hence, if $|I_s|<|H|$, we have
\[
 \dist(H,I_s)
 \ge\dist(z,H)-|I_s|
 \ge4|H|-|I_s|
 \ge3|H|,
\]
which proves \eqref{eq:whitney-separation}.
\end{proof}

\begin{lem}
\label{lem:truncation}
Let ${\bf T}$ be a tree that is lacunary in the first coordinate.
Assume that, for every $s\in{\bf T}$,
\[
 I_s\not\subset\Omega,\qquad  \textrm{and} \qquad |\beta_s|\le|I_s|.
\]
Then, for every $H\in\mathcal W$ with $|H|\le|I_{\bf T}|$, there holds
\begin{equation}\label{eq:local-truncation}
 \norm{\Psi_{\bf T}}_{L^\infty(H)}
 \lesssim
 1+\inf_{y\in H}M\Psi_{\bf T}(y).
\end{equation}
\end{lem}

\begin{proof}
Fix an $H\in\mathcal W$ with $h=|H|\le|I_{\bf T}|$, and choose $\xi_{{\bf T},1}$ satisfying the lacunarity
condition in Definition~\ref{def:tree}.
We split the packet sum into two parts
\[
 \Psi_{\bf T}
 =\Psi_{{\bf T},H}^{\ge}+\Psi_{{\bf T},H}^{<},
\]
where
\[
 \Psi_{{\bf T},H}^{\ge}
 =\sum_{\substack{s\in{\bf T}\\|I_s|\ge h}}
       \beta_s\varphi_{s,1},
 \qquad \textrm{and} \qquad 
 \Psi_{{\bf T},H}^{<}
 =\sum_{\substack{s\in{\bf T}\\|I_s|<h}}
       \beta_s\varphi_{s,1}.
\]

We first estimate $\Psi_{{\bf T},H}^{<}$ on $H$. Fix $\ell<h$ such that $|I_s|=\ell$ for some $s\in{\bf T}$.
By \eqref{eq:whitney-separation}, for every $s\in{\bf T}$
with $|I_s|=\ell$, we have $\dist(H,I_s)\ge3h$.
To estimate the contribution of these tri-tiles,
fix $x\in H$ and group them according to their distance
from $x$, namely
$n\ell\le\dist(x,I_s)<(n+1)\ell$, where $n\in\N$.
By Lemmas~\ref{lem:tree-intervals} and~\ref{lem:tree-counting},
each group contains only a bounded number of tri-tiles,
uniformly in $x$, $\ell$, and $n$.
Using \eqref{eq:packet} and $|\beta_s|\le\ell$, we therefore have
\[
 \sup_{x\in H}
 \sum_{\substack{s\in{\bf T}\\|I_s|=\ell}}
       |\beta_s\varphi_{s,1}(x)|
 \lesssim
 \sum_{\substack{n\in\N\\n\ge h/\ell}}(1+n)^{-A}
 \lesssim(\ell/h)^{A-1},
\]
where we recall $A \gg 1$ is a sufficiently large integer selected in \eqref{eq:packet}.  Summing over the spatial lengths $\ell<h$ and using Lemma~\ref{lem:tree-intervals},
we see that
\begin{equation}\label{eq:fine-physical}
 \norm{\Psi_{{\bf T},H}^{<}}_{L^\infty(H)}\lesssim1.
\end{equation}

To estimate $\Psi_{{\bf T},H}^{\ge}$, choose
$\eta\in C_c^\infty(\R)$ such that $\eta=1$ on $[-K,K]$
and $\supp\,\eta\subset[-2K,2K]$, where $K$ is the
constant in \eqref{eq:topfreq}.
Define the Fourier multiplier $P_{{\bf T},h}$ by
\[
\widehat{P_{{\bf T},h}u}(\xi):=
 \eta\bigl(h(\xi-\xi_{{\bf T},1})\bigr)\widehat u(\xi).
\]
For every $s\in{\bf T}$ with $|I_s|\ge h$,
\eqref{eq:topfreq} shows that the multiplier equals one
on $\supp\,\widehat\varphi_{s,1}$.
Consequently, one has
\begin{equation}\label{eq:cutoff-identity}
 \Psi_{{\bf T},H}^{\ge}
 = P_{{\bf T},h}\Psi_{\bf T}
 - P_{{\bf T},h}\Psi_{{\bf T},H}^{<}.
\end{equation}

We next estimate $P_{{\bf T},h}\Psi_{{\bf T},H}^{<}$. Set
\[
 J_h=
 \left[\xi_{{\bf T},1}-\frac{2K}{h},
       \xi_{{\bf T},1}+\frac{2K}{h}\right].
\]
We first show that only a bounded number of distinct
frequency intervals can contribute to
$P_{{\bf T},h}\Psi_{{\bf T},H}^{<}$.
Indeed, since the multiplier defining $P_{{\bf T},h}$
is supported in $J_h$ and
$\supp\,\widehat{\varphi}_{s,1}\subset\omega_{s,1}$,
we have $P_{{\bf T},h}\varphi_{s,1}=0$ whenever
$\omega_{s,1}\cap J_h=\varnothing$.
It therefore suffices to count the distinct intervals
$\omega_{s,1}$ with $s\in{\bf T}$ and $|I_s|<h$
that meet $J_h$.

For these tri-tiles, we have
$|\omega_{s,1}|=|I_s|^{-1}>h^{-1}$.
Moreover, by \eqref{eq:lacunary}, for any $s,s'\in{\bf T}$
with $\omega_{s,1}\ne\omega_{s',1}$, we have
$\omega_{s,1}\cap\omega_{s',1}=\varnothing$.
Since $|J_h|=4K/h$, we conclude that only
$\lesssim_K 1$ distinct frequency intervals can contribute.

Next, observe that the convolution kernel of $P_{{\bf T},h}$ is given by
\[
 K_{{\bf T},h}(z)
 :=
 h^{-1}e^{2\pi i\xi_{{\bf T},1}z}\check\eta(z/h).
\]
Choose $B>2$ with $B+2<A$, where $A \gg 1$ is defined as in \eqref{eq:packet}.
Then for every $s\in{\bf T}$ with $\ell=|I_s|<h$,
the decay of the packet and the kernel gives
\begin{align}
 |P_{{\bf T},h}\varphi_{s,1}(x)|
 &\lesssim
 \int_\R
 \frac1h\left(1+\frac{|x-z|}{h}\right)^{-B-2}
 \frac1\ell
 \left(1+\frac{|z-c(I_s)|}{\ell}\right)^{-A}\,dz
 \notag\\
 &\lesssim
 \frac1h
 \left(1+\frac{\dist(x,I_s)}h\right)^{-B}.
 \label{eq:smoothed-packet}
\end{align}
Here, in the last estimate above, we used the elementary estimate
\[
 1+\frac{|x-c(I_s)|}{h}
 \le
 \left(1+\frac{|x-z|}{h}\right)
 \left(1+\frac{|z-c(I_s)|}{\ell}\right).
\]

Fix one of the contributing frequency intervals $\omega$.
All corresponding spatial intervals have length
$\ell=|\omega|^{-1}<h$.
Thus, by \eqref{eq:smoothed-packet} and $|\beta_s|\le\ell$,
\begin{align*}
 \sum_{\substack{s\in{\bf T}\\\omega_{s,1}=\omega}}
       |\beta_sP_{{\bf T},h}\varphi_{s,1}(x)| \lesssim
 \sum_{\substack{s\in{\bf T}\\\omega_{s,1}=\omega}}
 \frac{\ell}{h}
 \left(1+\frac{\dist(x,I_s)}h\right)^{-B} \lesssim 1,
\end{align*}
where in the above estimate, we used the disjointness of the
distinct spatial intervals of length $\ell$ and their
bounded multiplicity from
Lemmas~\ref{lem:tree-intervals} and~\ref{lem:tree-counting}.
Summing over the bounded number of frequency intervals, we have
\begin{equation}\label{eq:fine-leakage}
 \norm{P_{{\bf T},h}\Psi_{{\bf T},H}^{<}}_\infty
 \lesssim1.
\end{equation}

Finally, we estimate $P_{{\bf T},h}\Psi_{\bf T}$.
For every $x,y\in H$, we have $|x-y|\le h$.
Using the decay of $K_{{\bf T},h}$ and decomposing
the integral into annuli centered at $y$, we find
\begin{align}
 |P_{{\bf T},h}\Psi_{\bf T}(x)|
 &\lesssim
 \frac1h\int_\R
 \left(1+\frac{|z-y|}{h}\right)^{-2}
 |\Psi_{\bf T}(z)|\,dz
 \notag\\
 &\lesssim
 \sum_{m\ge0}2^{-m}
 \frac1{2^{m+1}h}
 \int_{|z-y|\le2^m h}|\Psi_{\bf T}(z)|\,dz
 \notag\\
 &\lesssim M\Psi_{\bf T}(y).
 \label{eq:truncation-maximal}
\end{align}

Since
$\Psi_{\bf T}=\Psi_{{\bf T},H}^{\ge}+\Psi_{{\bf T},H}^{<}$,
we may apply \eqref{eq:cutoff-identity} to obtain, for every $x,y\in H$,
\begin{align*}
 |\Psi_{\bf T}(x)|
\le
 |P_{{\bf T},h}\Psi_{\bf T}(x)|
 +|P_{{\bf T},h}\Psi_{{\bf T},H}^{<}(x)|
 +|\Psi_{{\bf T},H}^{<}(x)|\lesssim 1+M\Psi_{\bf T}(y).
\end{align*}
For the last inequality, we applied
\eqref{eq:truncation-maximal}, \eqref{eq:fine-leakage},
and \eqref{eq:fine-physical},
respectively.
Taking the supremum over $x\in H$ and the infimum
over $y\in H$, we conclude \eqref{eq:local-truncation}.
\end{proof}

We now sum the local suprema over the Whitney intervals
for a single tree. 

\begin{lem}\label{lem:local-tree}
Let ${\bf T}$ be a tree such that
$I_s\not\subset\Omega$ for every $s\in{\bf T}$.
Suppose that
$\operatorname{size}_2(f_2;{\bf T})\le1$ and
$\operatorname{size}_3(f_3;{\bf T})\le1$.
Then
\begin{equation}\label{eq:local-tree}
 \sum_{H\in\mathcal W}|H|
       \norm{\Psi_{\bf T}}_{L^\infty(H)}
 \lesssim |I_{\bf T}|.
\end{equation}
\end{lem}

\begin{proof}
Write $L=|I_{\bf T}|$ and choose
$j\in\{1,2,3\}$ such that ${\bf T}$ is a $j$-tree.
We consider separately the cases for $j=1$ and $j\ne1$.

\smallskip
\noindent $\bullet$ \textsf{Case 1: $j=1$.}
By Lemma~\ref{lem:tree-frequencies}, the second and
third coordinates are lacunary.
Using Cauchy--Schwarz and \eqref{eq:size-coefficients},
\begin{equation}\label{eq:beta-l1}
 \begin{split}
 \sum_{s\in{\bf T}}|\beta_s|
 &\le
 \prod_{i=2}^3
 \left(
   \sum_{s\in{\bf T}}|I_s|\,
       |\langle f_i,\varphi_{s,i}\rangle|^2
 \right)^{1/2}\\
 &\le L.
 \end{split}
\end{equation}
Fix $s\in{\bf T}$, put $\ell=|I_s|$, and define
\[
 w_s(x)
 =
 \ell^{-1}
 \left(1+\frac{\dist(x,I_s)}{\ell}\right)^{-A}.
\]
If $H\in\mathcal W$ satisfies $|H|\le\ell$,
the values of $w_s$ on $H$ are comparable.
If $|H|>\ell$, the same conclusion follows from
$\dist(H,I_s)\ge3|H|$ in
\eqref{eq:whitney-separation}.
Since the Whitney intervals are disjoint,
\begin{align*}
 \sum_{H\in\mathcal W}|H|
       \norm{\varphi_{s,1}}_{L^\infty(H)}
 \lesssim
 \sum_{H\in\mathcal W}\int_H w_s(x)\,dx \le\int_\R w_s(x)\,dx
 \lesssim1.
\end{align*}
Multiplying by $|\beta_s|$, summing over $s\in{\bf T}$,
and using \eqref{eq:beta-l1}, we obtain
\eqref{eq:local-tree}.

\smallskip
\noindent $\bullet$ \textsf{Case 2: $j\ne1$.}
By Lemma~\ref{lem:tree-frequencies}, ${\bf T}$ is
lacunary in the first coordinate.
We first consider the Whitney intervals $H$ with
$|H|\le L$ and $H\cap8I_{\bf T}\ne\varnothing$.
To estimate their contribution, we combine
Lemma~\ref{lem:truncation} with an $L^2$ bound for
$\Psi_{\bf T}$. Applying Lemma~\ref{lem:forest2} to the forest
consisting of ${\bf T}$ alone, we have
\begin{equation}\label{eq:tree-l2}
 \norm{\Psi_{\bf T}}_2^2\lesssim L.
\end{equation}
For every $H\in\mathcal W$ with $|H|\le L$ and
$H\cap8I_{\bf T}\ne\varnothing$, we also have, by
Lemma~\ref{lem:truncation},
\[
 \norm{\Psi_{\bf T}}_{L^\infty(H)}
 \lesssim1+\inf_{y\in H}M\Psi_{\bf T}(y).
\]
Since these intervals are disjoint and contained in
$10I_{\bf T}$, we may sum over them and then apply
Cauchy--Schwarz and
\eqref{eq:tree-l2} to obtain
\begin{align}
 \sum_{\substack{H\in\mathcal W\\
          |H|\le L,\ H\cap8I_{\bf T}\ne\varnothing}}
       |H|\norm{\Psi_{\bf T}}_{L^\infty(H)}
 &\lesssim
 L+\int_{10I_{\bf T}}M\Psi_{\bf T}(x)\,dx
 \notag\\
 &\lesssim
 L+L^{1/2}\norm{\Psi_{\bf T}}_2
 \lesssim L.
 \label{eq:nearby-whitney}
\end{align}

We now turn to the remaining Whitney intervals, namely
\[
 \mathcal W_{\mathrm{far}}
 =
 \{H\in\mathcal W:
       |H|>L\ \text{or}\ H\cap8I_{\bf T}=\varnothing\}.
\]
To use packet decay on these intervals, we first estimate
their distance from $I_{\bf T}$.
For every $H\in\mathcal W_{\mathrm{far}}$, set
$D_H=\dist(H,I_{\bf T})$.
If $|H|\le L$, then $H\cap8I_{\bf T}=\varnothing$,
so $D_H\ge7L/2$.
If $|H|>L$, choose $s\in{\bf T}$ and
$z\in I_s\setminus\Omega$.
Since $z\in I_{\bf T}$ and $z\notin9H$, we have
\[
 D_H\ge\dist(z,H)-L
 \ge4|H|-L\ge3|H|.
\]
Thus, in both cases,
\begin{equation}\label{eq:far-whitney-distance}
 D_H\ge3\max(|H|,L).
\end{equation}

For every $H\in\mathcal W_{\mathrm{far}}$, we estimate
$\Psi_{\bf T}$ on $H$ by summing the packet bounds
over spatial lengths.
By Lemma~\ref{lem:tree-counting}, for each $\ell\le L$,
there are at most $CL/\ell$ tri-tiles $s\in{\bf T}$
with $|I_s|=\ell$.
Using this bound together with $|\beta_s|\le|I_s|$
and \eqref{eq:packet}, we have
\[
 \norm{\Psi_{\bf T}}_{L^\infty(H)}
 \lesssim
 \sum_{\ell\le L}\frac L\ell
       \left(\frac{\ell}{D_H}\right)^A
 \lesssim\left(\frac L{D_H}\right)^A.
\]
Here, in the above estimate, the sum is over the distinct spatial lengths in ${\bf T}$, and the last inequality follows from
Lemma~\ref{lem:tree-intervals}.

To sum over $H\in\mathcal W_{\mathrm{far}}$, note that by
\eqref{eq:far-whitney-distance}, we have for every
$x\in H$,
\[
 D_H\le\dist(x,I_{\bf T})
 \le D_H+|H|\le\frac43D_H.
\]
Since the Whitney intervals are disjoint and $A>1$,
we therefore conclude that
\begin{align*}
 \sum_{H\in\mathcal W_{\mathrm{far}}}
       |H|\norm{\Psi_{\bf T}}_{L^\infty(H)}
 &\lesssim
 \int_{\dist(x,I_{\bf T})\ge3L}
       \left(\frac L{\dist(x,I_{\bf T})}\right)^A\,dx\\
 &\lesssim L.
\end{align*}
Combining the estimate over $\mathcal W_{\mathrm{far}}$
with \eqref{eq:nearby-whitney}, we conclude
\eqref{eq:local-tree} when $j\ne1$.
\end{proof}

Finally, we complete the proof of Proposition~\ref{prop:local}.

\begin{proof}[Proof of Proposition~\ref{prop:local}]
Since size cannot increase under restriction, the
hypotheses of Lemma~\ref{lem:local-tree} hold on
every ${\bf T}\in\mathcal F$.
Using $\Psi=\sum_{{\bf T}\in\mathcal F}\Psi_{\bf T}$, we see that
\begin{align*}
 \sum_{H\in\mathcal W}|H|\norm{\Psi}_{L^\infty(H)}
\le
 \sum_{{\bf T}\in\mathcal F}
 \sum_{H\in\mathcal W}|H|
       \norm{\Psi_{\bf T}}_{L^\infty(H)} \lesssim
 \sum_{{\bf T}\in\mathcal F}|I_{\bf T}|. 
\end{align*}
The proof of Proposition~\ref{prop:local} is complete. 
\end{proof}

\subsection{Proof of the forest estimate}\label{subsec:forest-proof}

We now return to Proposition~\ref{prop:forest-p}.
Recall that
\[
 E_f=\{M_pf>C_0\},
 \qquad \textrm{and} \qquad 
 N=\sum_{{\bf T}\in\mathcal F}|I_{\bf T}|.
\]
Under the hypotheses of Proposition~\ref{prop:forest-p},
our goal is to prove
\[
 |\Lambda_{\bf S}(f,f_2,f_3)|
 \lesssim_{p,C_0}
 (N^{1/2}+N^{1/p})
 \operatorname{size}_2(f_2;{\bf S})
 \operatorname{size}_3(f_3;{\bf S}).
\]

\begin{proof}[Proof of Proposition~\ref{prop:forest-p}]
If either the second or third size is zero, every
summand in \eqref{eq:model} vanishes.
Otherwise, by homogeneity we may assume
\[
 \operatorname{size}_2(f_2;{\bf S})
 =
 \operatorname{size}_3(f_3;{\bf S})
 =1.
\]
Let $\Psi$ be the packet sum in \eqref{eq:synthesis}.
By \eqref{eq:forest-pairing}, it suffices to prove
\[
 |\langle f,\Psi\rangle|
 \lesssim_{p,C_0}(N^{1/2}+N^{1/p}).
\]
We estimate the pairing over $E_f^c$ first, and then
over the Whitney intervals of $E_f$.

On $E_f^c$, we have $|f|\le C_0$ almost everywhere,
and hence
\[
 \norm{f\one_{E_f^c}}_2^2
 \le C_0^{2-p}\norm{f}_p^p
 \le C_0^{2-p}.
\]
Applying \eqref{eq:forest-l2} with $h=f\one_{E_f^c}$,
we deduce that
\begin{equation}\label{eq:good-pairing}
 |\langle f\one_{E_f^c},\Psi\rangle|
 \lesssim_{p,C_0}N^{1/2}.
\end{equation}

We next consider the pairing over $E_f$.
The set $E_f$ is open, and by the weak $(1,1)$ bound
for $M$, we have
\[
 |E_f|
 \lesssim C_0^{-p}\norm{f}_p^p
 \le C_0^{-p}.
\]
Let $\mathcal W$ consist of the maximal intervals
$H$ of a fixed dyadic lattice such that $9H\subset E_f$,
as in Proposition~\ref{prop:local} with $\Omega=E_f$.
We first estimate the $L^p$ norm of $f$ on each
$H\in\mathcal W$.
By maximality, the dyadic parent $\widehat H$ satisfies
$9\widehat H\not\subset E_f$.
Choose $y\in9\widehat H\setminus E_f$.
Since $M_pf(y)\le C_0$ and $|9\widehat H|=18|H|$,
we have
\begin{equation}\label{eq:local-f}
 \int_H|f(x)|^p\,dx
 \le\int_{9\widehat H}|f(x)|^p\,dx
 \le18C_0^p|H|.
\end{equation}

To make use of \eqref{eq:local-f},
we estimate the $L^{p'}$ norm of $\Psi$ on each
Whitney interval.
Since $p'>2$, we have, for every $H\in\mathcal W$,
\begin{equation}\label{eq:local-interpolation}
 \norm{\Psi}_{L^{p'}(H)}
 \le
 \norm{\Psi}_{L^\infty(H)}^{2/p-1}
 \norm{\Psi}_{L^2(H)}^{2/p'}.
\end{equation}
Applying H\"older on each $H \in \mathcal W$, then using
\eqref{eq:local-f} and \eqref{eq:local-interpolation}, we have 
\begin{align}
 |\langle f\one_{E_f},\Psi\rangle|
 &\le
 \sum_{H\in\mathcal W}
       \norm{f}_{L^p(H)}\norm{\Psi}_{L^{p'}(H)}
 \notag\\
 &\lesssim_{p,C_0}
 \sum_{H\in\mathcal W}
       \left(|H|\norm{\Psi}_{L^\infty(H)}\right)^{2/p-1}
       \left(|H|^{1/2}\norm{\Psi}_{L^2(H)}\right)^{2/p'}
 \notag\\
 &\lesssim_{p,C_0}
 \left(
   \sum_{H\in\mathcal W}|H|\norm{\Psi}_{L^\infty(H)}
 \right)^{2/p-1}
 \left(
   \sum_{H\in\mathcal W}|H|^{1/2}\norm{\Psi}_{L^2(H)}
 \right)^{2/p'}.
 \label{eq:local-holder}
\end{align}
Here we used
$1/p=(2/p-1)+1/p'$ in the second step and
$(2/p-1)+2/p'=1$ in the last step.
If $\mathcal W$ is infinite, we may first estimate finite
subcollections and then pass to the limit.

For the first sum in \eqref{eq:local-holder},
we apply Proposition~\ref{prop:local} to obtain
\[
 \sum_{H\in\mathcal W}|H|\norm{\Psi}_{L^\infty(H)}
 \lesssim N.
\]
For the second sum, by the disjointness of the Whitney intervals, Cauchy--Schwarz, and \eqref{eq:coeff}, one has
\begin{align*}
 \sum_{H\in\mathcal W}|H|^{1/2}\norm{\Psi}_{L^2(H)}
 &\le
 |E_f|^{1/2}
 \left(
   \sum_{H\in\mathcal W}\norm{\Psi}_{L^2(H)}^2
 \right)^{1/2}\\
 &\le |E_f|^{1/2}\norm{\Psi}_2
 \lesssim_{p,C_0}N^{1/2},
\end{align*}
where we used $|E_f|\lesssim C_0^{-p}$.
Substituting the estimates for the two sums into
\eqref{eq:local-holder}, we conclude that
\begin{equation}\label{eq:bad-pairing}
 |\langle f\one_{E_f},\Psi\rangle|
 \lesssim_{p,C_0}
 N^{2/p-1}N^{1/p'}
 =N^{1/p}.
\end{equation}

Combining \eqref{eq:good-pairing} and
\eqref{eq:bad-pairing}, we deduce that
\[
 |\langle f,\Psi\rangle|
 \lesssim_{p,C_0}N^{1/2}+N^{1/p}.
\]
By \eqref{eq:forest-pairing}, we have therefore proved
\eqref{eq:forest-p} when the second and third sizes
equal one.
Multiplying by the original values of
$\operatorname{size}_2(f_2;{\bf S})$ and
$\operatorname{size}_3(f_3;{\bf S})$, we obtain  \eqref{eq:forest-p} in full generality.
\end{proof}

\bigskip 

\section{Proof of the model estimate}\label{sec:model-proof}

We now prove Proposition~\ref{prop:model} thereby completing the proof of Proposition~\ref{prop:restricted}.
Throughout this section, we assume \eqref{eq:normalized-inputs},
retain the sets $E,\widetilde E,F_3'$ from
\eqref{eq:exceptional}, and take
\[
 |f_2|\le\one_{F_2},\qquad |f_3|\le\one_{F_3'}.
\]
For a finite well-discretized collection ${\bf S}$, our goal
is to establish
\[
 |\Lambda_{\bf S}(f_1,f_2,f_3)|
 \lesssim_{p_1,p_2}|F_2|^{1/p_2}.
\]
We decompose ${\bf S}$ into
\[
 {\bf S}_{\mathrm{in}}
 =\{s\in{\bf S}:I_s\subset E\},
 \qquad
 {\bf S}_{\mathrm{out}}
 ={\bf S}\setminus{\bf S}_{\mathrm{in}},
\]
so that
\begin{equation}\label{eq:model-in-out}
 \Lambda_{\bf S}(f_1,f_2,f_3)
 =
 \Lambda_{{\bf S}_{\mathrm{in}}}(f_1,f_2,f_3)
 +
 \Lambda_{{\bf S}_{\mathrm{out}}}(f_1,f_2,f_3).
\end{equation}
For ${\bf S}_{\mathrm{in}}$, we use the vanishing of $f_3$
near the spatial intervals and the decay of the packets.
For ${\bf S}_{\mathrm{out}}$, we apply the forest estimate
from Proposition~\ref{prop:forest-p}.

\subsection{Tiles inside the exceptional set}\label{subsec:inside}

We follow the exceptional-set argument
in~\cite[Section~6, Proposition~6.1]{DPT}.

\begin{lem}\label{lem:interior}
Under the assumptions of Proposition~\ref{prop:model}, we have
\begin{equation}\label{eq:interior}
 \sum_{s\in{\bf S}_{\mathrm{in}}}
 |I_s|\prod_{i=1}^3|\langle f_i,\varphi_{s,i}\rangle|
 \lesssim_{p_1,p_2}|F_2|^{1/p_2}.
\end{equation}
\end{lem}

\begin{proof}
Let $Q$ range over the connected components of $E$.
Each $Q$ is a bounded interval with an endpoint outside $E$.
By \eqref{eq:exceptional}, for every $m\ge0$, we have
\begin{equation}\label{eq:interior-component-averages}
 \int_{2^mQ}|f_1|^{p_1}\lesssim2^m|Q|,
 \qquad
 \int_{2^mQ}|f_2|^{p_2}
 \lesssim |F_2|\,2^m|Q|.
\end{equation}
Moreover, $f_3$ vanishes on $3Q$ by
Lemma~\ref{lem:major-subset}.

Fix a spatial interval $I\subset Q$ such that $I=I_s$ for some $s\in{\bf S}_{\mathrm{in}}$, and set
$w_I(x)=(1+|x-c(I)|^2/|I|^2)^2$.
Since multiplication by the polynomial $w_I$ does not
enlarge Fourier support, we have
\[
 \supp\,\widehat{w_I\varphi_{s,i}}\subset\omega_{s,i}.
\]
Also, for each fixed $i$, the intervals $\omega_{s,i}$ with
$I_s=I$ are pairwise disjoint by
Definition~\ref{defn: well-discretized},   (iii).
Hence the packets $w_I\varphi_{s,i}$ are orthogonal. Moreover, we have
\[
 \norm{w_I\varphi_{s,i}}_\infty\lesssim |I|^{-1},
 \qquad
 \norm{w_I\varphi_{s,i}}_2\lesssim |I|^{-1/2}, 
\]
which, together with orthogonality, gives
\[
 \sup_{\substack{s\in{\bf S}_{\mathrm{in}}\\I_s=I}}
 |\langle f,w_I\varphi_{s,i}\rangle|
 \lesssim |I|^{-1}\norm{f}_1
\]
and
\[
 \left(
 \sum_{\substack{s\in{\bf S}_{\mathrm{in}}\\I_s=I}}
 |\langle f,w_I\varphi_{s,i}\rangle|^2
 \right)^{1/2}
 \lesssim |I|^{-1/2}\norm{f}_2.
\]
Interpolating these bounds and using
$\langle f,\varphi_{s,i}\rangle
=\langle f/w_I,w_I\varphi_{s,i}\rangle$,
we obtain, for $1<q\le2$,
\begin{equation}\label{eq:interior-coefficient}
 \left(
 \sum_{\substack{s\in{\bf S}_{\mathrm{in}}\\I_s=I}}
 |\langle f,\varphi_{s,i}\rangle|^{q'}
 \right)^{1/q'}
 \lesssim |I|^{-1/q}\norm{f/w_I}_q.
\end{equation}

Since $1/p_1'+1/p_2'+1/2=1$, we apply
H\"older, \eqref{eq:interior-coefficient}, and \eqref{eq:interior-component-averages}, to obtain
\begin{align}
 \sum_{\substack{s\in{\bf S}_{\mathrm{in}}\\I_s=I}}
 |I_s|\prod_{i=1}^3|\langle f_i,\varphi_{s,i}\rangle|
 &\lesssim
 |I|
 \prod_{i=1}^2
 \left(\frac1{|I|}\int_\R |f_i|^{p_i}w_I^{-p_i}\right)^{1/p_i}
 \left(\frac1{|I|}\int_{(3Q)^c}w_I^{-2}\right)^{1/2}
 \notag\\
 &\lesssim
 |I|
 \prod_{i=1}^2
 \left(
   \frac1{|I|}
   \sum_{m\ge0}2^{-4mp_i}
   \int_{2^{m+1}Q}|f_i|^{p_i}
 \right)^{1/p_i}
 \left(\frac{|I|}{|Q|}\right)^{7/2}
 \notag\\
 &\lesssim
 |F_2|^{1/p_2}|I|
 \left(\frac{|Q|}{|I|}\right)^{3/2}
 \left(\frac{|I|}{|Q|}\right)^{7/2} \cdot 
 \prod_{i=1}^2
 \left(\sum_{m\ge0}2^{-m(4p_i-1)}\right)^{1/p_i}
 \notag\\
 &\lesssim
 |F_2|^{1/p_2}|I|
 \left(\frac{|I|}{|Q|}\right)^2.
 \label{eq:interior-one-interval}
\end{align}
Here, for the integral over $(3Q)^c$, we used
$|x-c(I)|\ge|Q|$ whenever $x\notin3Q$.

Finally, by Definition~\ref{defn: well-discretized}, we sum
\eqref{eq:interior-one-interval} first over the spatial intervals
$I\subset Q$ and then over the components $Q$ to obtain
\[
 \sum_{s\in{\bf S}_{\mathrm{in}}}
 |I_s|\prod_{i=1}^3|\langle f_i,\varphi_{s,i}\rangle|
 \lesssim
 |F_2|^{1/p_2}\sum_Q|Q|
 =
 |F_2|^{1/p_2}|E|
 \lesssim |F_2|^{1/p_2},
\]
where the last inequality follows from
Lemma~\ref{lem:major-subset}. 
\end{proof}

\subsection{Tiles outside the exceptional set}\label{subsec:outside}

For the exterior tiles, we follow the argument in~\cite[Section~5.4]{DDP}
and~\cite[Section~7]{DPT}.
The key new ingredient is our proof of the forest estimate
in Proposition~\ref{prop:forest-p} using the local supremum
estimate. We apply Proposition~\ref{prop:forest-p} to each
selected forest.

\begin{lem}\label{lem:exterior}
Under the assumptions of Proposition~\ref{prop:model}, we have
\begin{equation}\label{eq:exterior}
 |\Lambda_{{\bf S}_{\mathrm{out}}}(f_1,f_2,f_3)|
 \lesssim_{p_1,p_2}|F_2|^{1/p_2}.
\end{equation}
\end{lem}

\begin{proof}
To prove \eqref{eq:exterior}, we will decompose
${\bf S}_{\mathrm{out}}$ into forests and apply
Proposition~\ref{prop:forest-p}. We divide the proof into three steps. 

\vspace{0.1cm}

\noindent \underline{\textit{Step I: Size bounds for the exterior collection ${\bf S}_{\mathrm{out}}$.}} 

\vspace{0.15cm}

 We first verify that the hypotheses of
Proposition~\ref{prop:forest-p} hold for every forest
selected from ${\bf S}_{\mathrm{out}}$.
Let $s\in{\bf S}_{\mathrm{out}}$, and hence the interval $I_s$ contains a point outside $E$.
By \eqref{eq:exceptional} and \eqref{eq:size-max},
we may fix a constant $C_1\ge1$, independent of the
input functions, such that
\begin{equation}\label{eq:exterior-sizes}
 \operatorname{size}_2(f_2;{\bf S}_{\mathrm{out}})
 \le C_1|F_2|,
 \qquad
 \operatorname{size}_3(f_3;{\bf S}_{\mathrm{out}})
 \le C_1.
\end{equation}
Moreover, for every $s\in{\bf S}_{\mathrm{out}}$,
\[
 I_s\not\subset\{M_{p_1}f_1>C\}.
\]
Thus, Proposition~\ref{prop:forest-p} applies to every
forest selected from ${\bf S}_{\mathrm{out}}$, with
$f=f_1$, $p=p_1$, and $C_0=C$.

\vspace{0.1cm}

\noindent \underline{\textit{Step II: Decomposition of the exterior form $\Lambda_{{\bf S}_{\mathrm{out}}}$.}}

\vspace{0.15cm}

Since $\norm{f_3}_2^2\le|F_3'|\le1$, we can use \eqref{eq:selection} in the third coordinate to control
the total top interval length at each size level.
Starting with ${\bf S}_{\mathrm{out}}$, we apply
\eqref{eq:selection} successively to the remaining
collection at thresholds $C_1 2^{-k}$, beginning with $k=0$.
For each $k\ge0$, let $\mathcal F_k$ be the forest
removed at threshold $C_1 2^{-k}$, and set
\[
 {\bf S}_k=\bigsqcup_{{\bf T}\in\mathcal F_k}{\bf T}.
\]
By \eqref{eq:selection}, we have
\begin{equation}\label{eq:k-levels}
 \operatorname{size}_3(f_3;{\bf S}_k)\le C_1 2^{-k},
 \qquad
 N_k:=\sum_{{\bf T}\in\mathcal F_k}|I_{\bf T}|
 \lesssim2^{2k}.
\end{equation}

The pairwise disjoint collections ${\bf S}_k$ cover
${\bf S}_{\mathrm{out}}$ except for a remainder which has zero
third size and hence has no contribution to $\Lambda_{{\bf S}_{\mathrm{out}}}$.
We therefore have
\begin{equation}\label{eq:exterior-level-decomposition}
 \Lambda_{{\bf S}_{\mathrm{out}}}(f_1,f_2,f_3)
 =
 \sum_{k\ge0}\Lambda_{{\bf S}_k}(f_1,f_2,f_3).
\end{equation}
Since ${\bf S}_{\mathrm{out}}$ is finite, only finitely
many collections ${\bf S}_k$ are nonempty, so the sum
in \eqref{eq:exterior-level-decomposition} is finite.

\vspace{0.1cm}

\noindent \underline{\textit{Step III: Summation over $k$.}}

\vspace{0.15cm}

Choose an integer $k_0\ge0$ such that
\begin{equation}\label{eq:k0}
 2^{-2k_0}\le|F_2|<4\cdot2^{-2k_0}.
\end{equation}
We estimate the ranges $0\le k\le k_0$ and $k>k_0$ in \eqref{eq:exterior-level-decomposition} separately.

\medskip 

\noindent $\bullet$ \textsf{Case 1: $0\le k\le k_0$.}
For each $0\le k\le k_0$, we apply
Proposition~\ref{prop:forest-p} to the forest $\mathcal F_k$.
Using \eqref{eq:exterior-sizes}, \eqref{eq:k-levels},
and the monotonicity of size, we have
\begin{equation}\label{eq:exterior-first-bound}
 |\Lambda_{{\bf S}_k}(f_1,f_2,f_3)|
 \lesssim_{p_1}
 |F_2| \cdot 2^{-k}(2^k+2^{2k/p_1})
 \lesssim_{p_1}
 |F_2| \cdot 2^{k(2/p_1-1)}.
\end{equation}
Then, summing
\eqref{eq:exterior-first-bound} over $0\le k\le k_0$,
we obtain
\begin{equation}\label{eq:small-k}
 \begin{split}
 \sum_{k=0}^{k_0}
       |\Lambda_{{\bf S}_k}(f_1,f_2,f_3)|
 &\lesssim_{p_1}
 |F_2| \cdot 2^{k_0(2/p_1-1)}\\
 &\lesssim_{p_1}
 |F_2|^{3/2-1/p_1}
 =
 |F_2|^{1/p_2}.
 \end{split}
\end{equation}
Here we used \eqref{eq:k0} for the second inequality
and \eqref{eq:critical} for the last equality.

\medskip 

\noindent $\bullet$ \textsf{Case 2: $k>k_0$.} In this case, we perform a size selection in the second coordinate.  Fix $k>k_0$.
We further decompose ${\bf S}_k$ in the second coordinate.
By \eqref{eq:exterior-sizes}, \eqref{eq:k0},
and the monotonicity of size, we have
\[
 \operatorname{size}_2(f_2;{\bf S}_k)
 \le C_1|F_2|
 \le2C_1 2^{-k_0}|F_2|^{1/2}.
\]
Starting with ${\bf S}_k$, we apply \eqref{eq:selection}
successively in the second coordinate at thresholds
\[
 2C_1 2^{-n}|F_2|^{1/2},
 \qquad n=k_0,\ldots,k-1.
\]
At each threshold, we apply the selection to the remaining collection.
For each $k_0\le n<k$, let $\mathcal F_{k,n}$ be
the forest removed at threshold
$2C_1 2^{-n}|F_2|^{1/2}$, and set
\[
 {\bf S}_{k,n}
 =
 \bigsqcup_{{\bf T}\in\mathcal F_{k,n}}{\bf T}.
\]
We denote the tri-tiles remaining after all these
selections by
\[
 {\bf S}_{k,k}
 :=
 {\bf S}_k\setminus
 \bigcup_{n=k_0}^{k-1}{\bf S}_{k,n}.
\]
Thus
\begin{equation}\label{eq:second-selection-decomposition}
 {\bf S}_k=\bigsqcup_{n=k_0}^{k}{\bf S}_{k,n}.
\end{equation}

Guided by \eqref{eq:second-selection-decomposition}, 
for every $k>k_0$, we split the sum to be estimated as
\begin{align}
 \sum_{k>k_0}
       |\Lambda_{{\bf S}_k}(f_1,f_2,f_3)|
 &\le
 \sum_{k>k_0}\sum_{n=k_0}^{k-1}
       |\Lambda_{{\bf S}_{k,n}}(f_1,f_2,f_3)|
 \notag\\
 &\quad+
 \sum_{k>k_0}
       |\Lambda_{{\bf S}_{k,k}}(f_1,f_2,f_3)|.
 \label{eq:large-k-decomposition}
\end{align}
For the first sum, each ${\bf S}_{k,n}$ already has
the forest decomposition $\mathcal F_{k,n}$.
To apply Proposition~\ref{prop:forest-p} to the second
sum, we also need a forest decomposition of
${\bf S}_{k,k}$.

Since ${\bf S}_{k,k}\subset{\bf S}_k$ and
${\bf S}_k$ is decomposed into the forest $\mathcal F_k$,
we have
$$
 {\bf S}_{k,k}
 =
 \bigsqcup_{{\bf T}\in\mathcal F_k}
       ({\bf T}\cap{\bf S}_{k,k}).
$$
For every nonempty intersection
${\bf T}\cap{\bf S}_{k,k}$, we retain the top interval
and frequency of ${\bf T}$.
By Definition~\ref{def:tree}, each ${\bf T} \cap {\bf S}_{k, k}$ is a tree, and hence
\[
 \mathcal F_{k,k}
 :=
 \left\{
 {\bf T}\cap{\bf S}_{k,k}:
 {\bf T}\in\mathcal F_k,\ 
 {\bf T}\cap{\bf S}_{k,k}\ne\varnothing
 \right\}
\]
is a forest decomposition of ${\bf S}_{k,k}$.

We can now apply Proposition~\ref{prop:forest-p}
to every term on the right side of
\eqref{eq:large-k-decomposition}.
For $k_0\le n\le k$, write
\[
 N_{k,n}
 :=
 \sum_{{\bf T}\in\mathcal F_{k,n}}|I_{\bf T}|.
\]
Returning to \eqref{eq:large-k-decomposition},
we obtain
\begin{align}
 \sum_{k>k_0}
       |\Lambda_{{\bf S}_k}(f_1,f_2,f_3)|
\lesssim_{p_1}
 \sum_{k>k_0}\sum_{n=k_0}^{k}
\bigl(N_{k,n}^{1/2}+N_{k,n}^{1/p_1}\bigr) \cdot 
 \prod_{i=2}^3
       \operatorname{size}_i(f_i;{\bf S}_{k,n}).
 \label{eq:large-k-forest-bound}
\end{align}

To estimate the double sum in
\eqref{eq:large-k-forest-bound}, we claim that,
for every $k_0\le n\le k$,
\begin{equation}\label{eq:kn-levels}
 \begin{split}
 \operatorname{size}_2(f_2;{\bf S}_{k,n})
 &\le 2C_1 2^{-n}|F_2|^{1/2},\\
 \operatorname{size}_3(f_3;{\bf S}_{k,n})
 \le C_1 2^{-k}, &
 \qquad\text{and}\qquad
 N_{k,n}\lesssim 2^{2n}.
 \end{split}
\end{equation}

\smallskip
\noindent\textit{Proof of the claim.}
For the first estimate, we consider first the case that $n<k$.
By the selection at threshold
$2C_1 2^{-n}|F_2|^{1/2}$ and the monotonicity of size,
we have
\[
 \operatorname{size}_2(f_2;{\bf S}_{k,n})
 \le2C_1 2^{-n}|F_2|^{1/2}.
\]
For $n=k$, we apply \eqref{eq:selection} at the last
threshold, corresponding to $n=k-1$, to obtain
\[
 \operatorname{size}_2(f_2;{\bf S}_{k,k})
 \le\frac12\left(2C_1 2^{-(k-1)}|F_2|^{1/2}\right)
 =2C_1 2^{-k}|F_2|^{1/2}.
\]

Next, we prove the second estimate.
For every $k_0\le n\le k$, we have
${\bf S}_{k,n}\subset{\bf S}_k$.
Hence, by the monotonicity of size and \eqref{eq:k-levels}, we have
\[
 \operatorname{size}_3(f_3;{\bf S}_{k,n})
 \le\operatorname{size}_3(f_3;{\bf S}_k)
 \le C_1 2^{-k}.
\]

Finally, we prove the last estimate.
For $n<k$, by \eqref{eq:selection} and
$\norm{f_2}_2^2\le|F_2|$, we see that 
\[
 N_{k,n}
 \lesssim
 \left(2C_1 2^{-n}|F_2|^{1/2}\right)^{-2}
 \norm{f_2}_2^2
 \lesssim2^{2n}.
\]
For $n=k$, note that each tree in $\mathcal F_{k,k}$ retains
the top interval of a distinct tree in $\mathcal F_k$.
Thus, by \eqref{eq:k-levels}, we have
\[
 N_{k,k}
 \le
 \sum_{{\bf T}\in\mathcal F_k}|I_{\bf T}|
 =N_k
 \lesssim2^{2k}.
\]
The proof of \eqref{eq:kn-levels} is complete. 

\medskip 

Substituting \eqref{eq:kn-levels} into
\eqref{eq:large-k-forest-bound} and summing first
over $n$, we find that
\begin{align}
 \sum_{k>k_0}
       |\Lambda_{{\bf S}_k}(f_1,f_2,f_3)|
 &\lesssim_{p_1}
 |F_2|^{1/2}
 \sum_{k>k_0}\sum_{n=k_0}^{k}
       2^{-n-k}(2^n+2^{2n/p_1})
 \notag\\
 &\lesssim_{p_1}
 |F_2|^{1/2}
 \sum_{k>k_0}2^{-k}
 \sum_{n=k_0}^{k}2^{n(2/p_1-1)}
 \notag\\
 &\lesssim_{p_1}
 |F_2|^{1/2}\sum_{k>k_0}2^{-2k/p_1'}
 \notag\\
 &\lesssim_{p_1}
 |F_2|^{1/2}2^{-2k_0/p_1'}
 \lesssim_{p_1}
 |F_2|^{1/2+1/p_1'}
 =
 |F_2|^{1/p_2}.
 \label{eq:large-k}
\end{align}
Here, the last inequality follows from \eqref{eq:k0}, and the last equality follows from \eqref{eq:critical}.

\vspace{0.1cm}

Finally, combining \eqref{eq:small-k} and \eqref{eq:large-k}
with \eqref{eq:exterior-level-decomposition},
we conclude \eqref{eq:exterior}.
\end{proof}

\subsection{Putting everything together}

We now combine the interior and exterior estimates
to prove Proposition~\ref{prop:model}.

\begin{proof}[Proof of Proposition~\ref{prop:model}]
By \eqref{eq:model-in-out} and
Lemmas~\ref{lem:interior} and~\ref{lem:exterior}, we have
\begin{align*}
 |\Lambda_{\bf S}(f_1,f_2,f_3)|
 \le
 |\Lambda_{{\bf S}_{\mathrm{in}}}(f_1,f_2,f_3)|
 +
 |\Lambda_{{\bf S}_{\mathrm{out}}}(f_1,f_2,f_3)|\lesssim_{p_1,p_2}|F_2|^{1/p_2}.
\end{align*}
This proves \eqref{eq:model-target}.
\end{proof}

Finally, combining Proposition~\ref{prop:model} with the
reductions in Sections~\ref{subsec:testing},
\ref{subsec:normalization}, and~\ref{subsec:model},
we complete the proof of Proposition~\ref{prop:restricted}.

\bigskip

\section{Completion of the proof}\label{sec:completion}

Finally, we prove  Theorem~\ref{thm:main}. 

\subsection{Weak-space summation and decomposition of the input functions}

The following result is due to
Stein, Taibleson, and Weiss~\cite{STW}, and for the reader's convenience, we include a short proof of the form used below.

\begin{lem}\label{lem:superposition}
Let $0<r<1$ and let $h_1,\ldots,h_J\in L^{r,\infty}(\R)$.
Then
\begin{equation}\label{eq:superposition}
 \left\|\sum_{j=1}^J h_j\right\|_{L^{r,\infty}}^r
 \le
 \frac{2-r}{1-r}
 \sum_{j=1}^J\norm{h_j}_{L^{r,\infty}}^r.
\end{equation}
\end{lem}

\begin{proof}
Fix $\lambda>0$. Note that
\[
 \left\{\left|\sum_{j=1}^J h_j\right|>\lambda\right\}
 \subset
 \bigcup_{j=1}^J\{|h_j|>\lambda\}
 \;\cup\;
 \left\{
 \left|\sum_{j=1}^J h_j\one_{\{|h_j|\le\lambda\}}\right|
 >\lambda
 \right\}.
\]
Therefore, by Chebyshev, we have 
\begin{align*}
 \left|\left\{\left|\sum_{j=1}^J h_j\right|>\lambda\right\}\right|
 &\le
 \sum_{j=1}^J|\{|h_j|>\lambda\}|
 +
 \left|\left\{
 \sum_{j=1}^J|h_j|\one_{\{|h_j|\le\lambda\}}>\lambda
 \right\}\right|\\
 &\le
 \sum_{j=1}^J|\{|h_j|>\lambda\}|
 +
 \frac1\lambda\sum_{j=1}^J
 \int_{\{|h_j|\le\lambda\}}|h_j(x)|\,dx\\
 &\le
 \sum_{j=1}^J|\{|h_j|>\lambda\}|
 +
 \frac1\lambda\sum_{j=1}^J
 \int_0^\lambda|\{|h_j|>t\}|\,dt\\
 &\le
 \lambda^{-r}\sum_{j=1}^J\norm{h_j}_{L^{r,\infty}}^r
 +
 \frac{\lambda^{-r}}{1-r}
 \sum_{j=1}^J\norm{h_j}_{L^{r,\infty}}^r,
\end{align*}
which clearly gives \eqref{eq:superposition}.
\end{proof}

We also record the following standard lemma for decomposing
each input into a sum of bounded functions supported on
sets of finite measure.

\begin{lem}\label{lem:quantiles}
Let $1\le p<\infty$ and $f\in L^p(\R)$.
There exist pairwise disjoint measurable sets $E_j$,
$j\in\Z$, and numbers $a_j\ge0$ such that, with
$f^{(j)}=f\one_{E_j}$,
$$
 f=\sum_{j\in\Z}f^{(j)},\qquad
 |E_j|\le2^j,\qquad
 |f^{(j)}|\le a_j\one_{E_j},
$$
and
\begin{equation}\label{eq:quantile-sum}
 \sum_{j\in\Z}2^ja_j^p\le4\norm{f}_p^p.
\end{equation}
The partial sums $\sum_{|j|\le N}f^{(j)}$ converge to $f$
almost everywhere and in $L^p(\R)$ as $N\to\infty$.
\end{lem}

\begin{proof}
For $t>0$, let
\[
 f^*(t)
 =
 \inf\{\lambda\ge0:|\{|f|>\lambda\}|\le t\}
\]
be the nonincreasing rearrangement of $|f|$.
For each $j\in\Z$, choose a measurable set $G_j$ such that
\[
 \{|f|>f^*(2^j)\}
 \subset G_j
 \subset \{|f|\ge f^*(2^j)\}\cap\{f\ne0\},
 \qquad \textrm{and} \qquad 
 |G_j|=\min\{2^j,|\{f\ne0\}|\}.
\]
Moreover, we may arrange that
$G_j\subset G_{j+1}$.

Set
\[
 E_j=G_j\setminus G_{j-1},
 \qquad
 a_j=f^*(2^{j-1}),
 \qquad
 f^{(j)}=f\one_{E_j}.
\]
The sets $E_j$ are pairwise disjoint, and
\[
 |E_j|=|G_j|-|G_{j-1}|\le2^{j-1}\le2^j.
\]
Since $G_{j-1}$ contains $\{|f|>f^*(2^{j-1})\}$,
we also have $|f^{(j)}|\le a_j\one_{E_j}$. This proves the bounds for $|E_j|$ and $|f^{(j)}|$.

Next, by the monotonicity of $f^*$, we obtain
\begin{align*}
 \sum_{j\in\Z}2^ja_j^p \le
 4\sum_{j\in\Z}
 \int_{2^{j-2}}^{2^{j-1}}f^*(t)^p\,dt =4\int_0^\infty f^*(t)^p\,dt
 =4\norm{f}_p^p,
\end{align*}
which proves \eqref{eq:quantile-sum}.

Finally, since $f^*(2^j)\to0$ as $j\to\infty$,
$\bigcup_{j\in\Z}G_j=\{f\ne0\}$ up to a null set.
Also, $|G_j|\le2^j\to0$ as $j\to-\infty$, and hence
$\bigcap_{j\in\Z}G_j$ has measure zero.
Therefore, the sets $E_j$ partition $\{f\ne0\}$ up to
a null set, and
\[
 \sum_{|j|\le N}f^{(j)}\longrightarrow f
 \qquad\text{a.e}.
\]
Since the sets $E_j$ are disjoint, we have
\[
 \left|f-\sum_{|j|\le N}f^{(j)}\right|^p
 \le |f|^p.
\]
Applying dominated convergence, we conclude that the
partial sums also converge to $f$ in $L^p(\R)$.
\end{proof}

\subsection{From the restricted estimate to Theorem~\ref{thm:main}}

\begin{proof}[Proof of Theorem~\ref{thm:main}]
Fix $p_1,p_2$ satisfying \eqref{eq:critical},
and let $f_1,f_2\in\calS(\R)$.
We divide the proof into several steps.

\vspace{0.1cm}

\noindent \underline{\textit{Step I: Decomposition of the input functions.}}

\vspace{0.15cm}

For each $i\in\{1,2\}$, apply
Lemma~\ref{lem:quantiles} to $f_i$ with exponent $p_i$.
Denote the resulting sets and amplitudes by
$E_{i,j}$ and $a_{i,j}$, and set
$f_i^{(j)}=f_i\one_{E_{i,j}}$.
Then
\begin{equation}\label{eq: inputbounds}
 f_i=\sum_{j\in\Z}f_i^{(j)},\qquad
 |f_i^{(j)}|\le a_{i,j}\one_{E_{i,j}},\qquad
 |E_{i,j}|\le2^j.
\end{equation}
By \eqref{eq:quantile-sum}, we have
\begin{equation}\label{eq:input-quantile-bounds}
 \sum_{j\in\Z}2^ja_{1,j}^{p_1}
 \le4\norm{f_1}_{p_1}^{p_1},
 \qquad
 \sum_{k\in\Z}2^ka_{2,k}^{p_2}
 \le4\norm{f_2}_{p_2}^{p_2}.
\end{equation}

For every integer $N\ge1$, by Lemma~\ref{lem:superposition},
we see that
\begin{equation}\label{eq:output-piece-sum}
 \left\|
 \textnormal{\bf BHT}
 \left(\sum_{|j|\le N}f_1^{(j)},
       \sum_{|k|\le N}f_2^{(k)}\right)
 \right\|_{L^{2/3,\infty}}^{2/3}
 \lesssim
 \sum_{\substack{|j|\le N\\|k|\le N}}
 \norm{\textnormal{\bf BHT}
       (f_1^{(j)},f_2^{(k)})}_{L^{2/3,\infty}}^{2/3}.
\end{equation}

\vspace{0.1cm}

\noindent \underline{\textit{Step II: Estimates of the right--hand side of \eqref{eq:output-piece-sum}.}}

\vspace{0.15cm}

Choose $\eta>0$ sufficiently small such that all four
numbers $1/p_1\pm\eta$ and $1/p_2\pm\eta$ lie strictly
between $1/2$ and $1$.
Note that
\begin{equation}\label{eq:nearby-exponents}
 \left(\frac1{p_1}+\eta\right)
 +\left(\frac1{p_2}-\eta\right)
 =
 \left(\frac1{p_1}-\eta\right)
 +\left(\frac1{p_2}+\eta\right)
 =\frac32.
\end{equation}
Thus, by Proposition~\ref{prop:restricted} and
\eqref{eq: inputbounds}, we have, for every $j,k\in\Z$,
\begin{equation}\label{eq:two-restricted}
 \begin{split}
 &\norm{\textnormal{\bf BHT}
       (f_1^{(j)},f_2^{(k)})}_{L^{2/3,\infty}}\\
 &\quad\lesssim_{p_1,p_2}a_{1,j}a_{2,k}
 \min\bigl\{
 |E_{1,j}|^{1/p_1+\eta}|E_{2,k}|^{1/p_2-\eta},
 |E_{1,j}|^{1/p_1-\eta}|E_{2,k}|^{1/p_2+\eta}
 \bigr\}.
 \end{split}
\end{equation}
Since $|E_{1,j}|\le2^j$ and $|E_{2,k}|\le2^k$,
we further have that
\begin{align}
 \norm{\textnormal{\bf BHT}
       (f_1^{(j)},f_2^{(k)})}_{L^{2/3,\infty}}
 &\lesssim_{p_1,p_2}
 (2^{j/p_1}a_{1,j})(2^{k/p_2}a_{2,k})
 \min\{2^{\eta(j-k)},2^{-\eta(j-k)}\}
 \notag\\
 &=
 2^{-\eta|j-k|}
 (2^{j/p_1}a_{1,j})(2^{k/p_2}a_{2,k}).
 \label{eq:off-diagonal}
\end{align}

\vspace{0.1cm}

\noindent \underline{\textit{Step III: Summation over the input pieces.}}

\vspace{0.15cm}

Substituting \eqref{eq:off-diagonal} into
\eqref{eq:output-piece-sum}, we have
\begin{equation}\label{eq:weighted-piece-sum}
 \begin{split}
 &\left\|
 \textnormal{\bf BHT}
 \left(\sum_{|j|\le N}f_1^{(j)},
       \sum_{|k|\le N}f_2^{(k)}\right)
 \right\|_{L^{2/3,\infty}}^{2/3}\\
 &\qquad\lesssim_{p_1,p_2}
 \sum_{\substack{|j|\le N\\|k|\le N}}
 2^{-2\eta|j-k|/3}
 (2^{j/p_1}a_{1,j})^{2/3}
 (2^{k/p_2}a_{2,k})^{2/3}.
 \end{split}
\end{equation}
To estimate the double sum, we first apply the change of variable $m=j-k$. Then, for every $m\in\Z$, by H\"older, we deduce that
\begin{equation}\label{eq:diagonal-holder}
 \begin{split}
 \sum_{\substack{|k|\le N\\|k+m|\le N}}
 (2^{(k+m)/p_1}a_{1,k+m})^{2/3}
 (2^{k/p_2}a_{2,k})^{2/3}
 \le
 \left(\sum_{|j|\le N}2^ja_{1,j}^{p_1}\right)^{2/(3p_1)}
 \left(\sum_{|k|\le N}2^ka_{2,k}^{p_2}\right)^{2/(3p_2)}.
 \end{split}
\end{equation}
Multiplying \eqref{eq:diagonal-holder} by
$2^{-2\eta|m|/3}$ and summing over $m\in\Z$, we obtain
\begin{align*}
 \textrm{RHS of \eqref{eq:weighted-piece-sum}}
 &=\sum_{\substack{|j|\le N\\|k|\le N}}
 2^{-2\eta|j-k|/3}
 (2^{j/p_1}a_{1,j})^{2/3}
 (2^{k/p_2}a_{2,k})^{2/3}\\
 &\le
 \left(\sum_{m\in\Z}2^{-2\eta|m|/3}\right)
 \left(\sum_{|j|\le N}2^ja_{1,j}^{p_1}\right)^{2/(3p_1)}
 \left(\sum_{|k|\le N}2^ka_{2,k}^{p_2}\right)^{2/(3p_2)}.
\end{align*}
Since $\sum_{m\in\Z}2^{-2\eta|m|/3}<\infty$,
we may substitute the above estimate into
\eqref{eq:weighted-piece-sum} and raise both sides
to the power $3/2$.
We conclude that
\begin{equation}\label{eq:finite-quantile-sum}
 \begin{split}
 &\left\|
 \textnormal{\bf BHT}
 \left(\sum_{|j|\le N}f_1^{(j)},
       \sum_{|k|\le N}f_2^{(k)}\right)
 \right\|_{L^{2/3,\infty}}\\
 &\qquad\lesssim_{p_1,p_2}
 \left(\sum_{|j|\le N}2^ja_{1,j}^{p_1}\right)^{1/p_1}
 \left(\sum_{|k|\le N}2^ka_{2,k}^{p_2}\right)^{1/p_2}.
 \end{split}
\end{equation}
By \eqref{eq:input-quantile-bounds}, the right side
of \eqref{eq:finite-quantile-sum} is bounded by
$C_{p_1,p_2}\norm{f_1}_{p_1}\norm{f_2}_{p_2}$,
independently of $N$.

\vspace{0.1cm}

\noindent \underline{\textit{Step IV: Passage to the limit and extension.}}

\vspace{0.15cm}

For each $i\in\{1,2\}$, the partial sums
$\sum_{|j|\le N}f_i^{(j)}$ converge to $f_i$
almost everywhere and are bounded in absolute value
by $|f_i|$.
Since $f_i\in L^2(\R)$, we may apply dominated
convergence to conclude that
\[
 \sum_{|j|\le N}f_i^{(j)}
 \longrightarrow f_i
 \qquad\text{in }L^2(\R).
\]
By bilinearity and the $L^2\times L^2\to L^1$
estimate for $\textnormal{\bf BHT}$, we have
\begin{align*}
 &\left\|
 \textnormal{\bf BHT}
 \left(\sum_{|j|\le N}f_1^{(j)},
       \sum_{|k|\le N}f_2^{(k)}\right)
 -
 \textnormal{\bf BHT}(f_1,f_2)
 \right\|_1\\
 &\qquad\lesssim
 \left\|\sum_{|j|\le N}f_1^{(j)}-f_1\right\|_2
 \left\|\sum_{|k|\le N}f_2^{(k)}\right\|_2+
 \norm{f_1}_2
 \left\|\sum_{|k|\le N}f_2^{(k)}-f_2\right\|_2
 \longrightarrow0.
\end{align*}
We may therefore pass to a subsequence along which
the outputs converge almost everywhere to
$\textnormal{\bf BHT}(f_1,f_2)$.
Applying Fatou's lemma and using
\eqref{eq:finite-quantile-sum} and
\eqref{eq:input-quantile-bounds}, we deduce that
\[
 \norm{\textnormal{\bf BHT}(f_1,f_2)}_{L^{2/3,\infty}}
 \lesssim_{p_1,p_2}
 \norm{f_1}_{p_1}\norm{f_2}_{p_2}.
\]
Thus \eqref{eq:main} holds for Schwartz inputs.
Finally, by a standard density argument, we extend the estimate
to all $f_1\in L^{p_1}(\R)$ and $f_2\in L^{p_2}(\R)$,
with $p_1,p_2$ satisfying \eqref{eq:critical}.
This completes the proof of Theorem~\ref{thm:main}.
\end{proof}

\bigskip 

\noindent{\bf Acknowledgement.}
The author was supported by the NSF grant DMS-2555999 and by the Simons Travel grant MPS-TSM-00007213.

\end{document}